\documentclass{amsart}
\usepackage[T1]{fontenc}
\usepackage[latin1]{inputenc}
\usepackage{latexsym,amssymb,amsmath}
\usepackage[pdftex,colorlinks=true]{hyperref}
\usepackage{amsthm}
\usepackage{longtable}
\usepackage{epsfig}
\usepackage{hhline}
\usepackage{epic}
\usepackage{times}

\newcommand{\Galg}{\mathbf{G}}
\newcommand{\Balg}{\mathbf{B}}
\newcommand{\Zalg}{\mathbf{Z}}
\newcommand{\Talg}{\mathbf{T}}

\newcommand{\PGL}{\operatorname{PGL}}

\newcommand{\F}{\mathbb{F}}
\newcommand{\N}{\mathbb{N}}
\newcommand{\Q}{\mathbb{Q}}
\newcommand{\Z}{\mathbb{Z}}
\newcommand{\cL}{\mathcal{L}}
\newcommand{\Hom}{\operatorname{Hom}}
\newcommand{\SL}{\operatorname{SL}}
\newcommand{\GL}{\operatorname{GL}}
\newcommand{\Id}{\operatorname{Id}}
\newcommand{\tr}{\textrm{tr}}
\newcommand{\lcm}{\operatorname{lcm}}

\newtheorem{theorem}{Theorem}[section] 
\newtheorem{lemma}[theorem]{Lemma}

\newtheorem{proposition}[theorem]{Proposition}

\newtheorem{definition}[theorem]{Definition}
\newtheorem{remark}[theorem]{Remark}

\title{On defining characteristic representations of finite reductive groups}

\author{Olivier Brunat}
\address{OB: 
UFR de Mathématiques\\
Université Denis Diderot - Paris 7\\
175, rue du Chevaleret\\
F-75013 Paris\\}
\email{brunat@math.jussieu.fr}

\author{Frank Lübeck}
\address{FL:
Lehrstuhl D für Mathematik\\
RWTH Aachen\\
Templergraben 64\\
D-52062 Aachen\\}
\email{Frank.Luebeck@Math.RWTH-Aachen.De}

\subjclass{20C15,\ 20C33}

\begin{document}

\begin{abstract} 
We  give  parameterizations of  the  irreducible  representations of  finite
groups of Lie type in their defining characteristic.
\end{abstract}

\maketitle

\section{Introduction}\label{intro}

We consider  series of finite  groups of Lie type  which are specified  by a
root datum  and a  finite order  automorphism of that  root datum.  For each
power $q$  of a prime  $p$ this  determines a connected  reductive algebraic
group $\Galg$  over $\bar\F_p$ (an  algebraic closure of the field  with $p$
elements) and a group of fixed  points $\Galg^F$ of a Frobenius morphism $F:
\Galg \to \Galg$, up to isomorphism.

We  are interested  in  a  parameterization of  the  irreducible modules  of
$\Galg^F$ over $\bar\F_p$.

A  known solution  to this  task is  to use  that the  groups $\Galg^F$  are
groups  with a  split $(B,N)$-pair  of  characteristic $p$.  There exists  a
parameterization of  the irreducible  representations over  an algebraically
closed field of  characteristic $p$ of such groups. For  details we refer to
the description by Curtis in~\cite[B,  Thm.5.7]{Bo70}. In our setup it would
be very technical  to construct the data for this  parameterization from the
given root  datum with  Frobenius action.  That parameterization  looks very
different for different Frobenius actions on the same algebraic group.

In the literature on representations of connected reductive algebraic groups
and finite groups  of Lie type in their defining characteristic most authors
restrict their descriptions to the case of simply connected algebraic groups
and the finite groups of Lie type arising from these.

In  this  case  there  is  a  nice  combinatorial  parameterization  of  the
absolutely irreducible modules of the algebraic group by the set of dominant
weights. The  irreducible modules of  the finite groups are  restrictions of
those  of the  algebraic group  and Steinberg~\cite{St68}  described a  nice
subset of dominant  weights which yields representatives  of the isomorphism
classes of  these modules.  A generalization  to connected  reductive groups
with simply connected derived group can be found in~\cite[App.1.3]{Hz09}.

Jantzen  considers  in~\cite{Ja03}  general  connected  reductive  algebraic
groups but does not  consider the finite groups of Lie  type. In the general
case it is no longer true that all irreducible representations of the finite
groups are restrictions from the algebraic group.

In this paper we give  a parameterization of the irreducible representations
in defining  characteristic for arbitrary finite  groups of Lie type.  It is
very concrete  and computable  starting from  the given  root datum  for the
algebraic group and Frobenius action on the root datum. The description will
not become more complicated for twisted Frobenius actions.

Here is  an overview  of the content  of the other  sections of  this paper.
Section~\ref{sec:rootdata} contains a description  of our setup. We describe
how root data and Frobenius actions on  root data can be represented and how
to compute certain related data. Some of  the results in this section may be
of  independent interest.  For  example,  we describe  a  construction of  a
certain  covering group  of an  arbitrary connected  reductive group,  which
generalizes the well  known simply connected coverings  of semisimple groups
(see~\ref{prop:covering}).

In section~\ref{sec:irrdefchar}  we first recall the  results about defining
characteristic representations of the algebraic groups and the finite groups
of Lie  type arising from simply  connected semisimple groups which  we have
mentioned above. Then we state  our main result in theorem~\ref{maintheorem}
where we consider  arbitrary finite groups of  Lie type. In the  end of that
section  we work  out an  example in  some detail  (certain centralizers  of
semisimple elements in exceptional groups of type $E_8$).

In  section~\ref{sec:simple} we  give  a more  detailed  description of  the
parameter sets  in our main  theorem for finite  groups of Lie  type arising
from  any simple  connected  reductive  group. As  an  application of  these
results we  work out the number  of semisimple conjugacy classes  for all of
these finite groups. The results of this application were obtained before by
the first  named author  with a  completely different  proof. The  new proof
given here is more elementary.

\textbf{Acknowledgements.}  We would  like to  thank Bob  Guralnick for  the
suggestion to combine a reduction  to the simply-connected case and Clifford
theory, as we  do in our main theorem~\ref{maintheorem}. We  also thank Marc
Cabanes for pointing  us to his result~\cite[B.11.3]{Cab88},  we have reused
his proof for our proposition~\ref{prop:ext}. In  the last section we need a
combinatorial lemma~\ref{lemma:an}, we thank Darij Grinberg for showing us a
proof, and for allowing us to include it in this paper.

\section{Root data for finite groups of Lie type}
\label{sec:rootdata}

\subsection{Connected reductive algebraic groups}
\label{ssec:rootdata}
Let  $\Galg$ be  a connected  reductive group  over an  algebraically closed
field $\bar  k$. We recall  how $\Galg$ is determined  by a root  datum, for
more details we refer to \cite[7.4, 9.6]{Sp98}.

For each maximal torus $\Talg$ of  $\Galg$ there is an associated root datum
$\Psi = (X,  R, Y, R^\vee)$ which together with  $\bar k$ determines $\Galg$
up  to isomorphism.  Here, $X$  and $Y$  are the  character and  cocharacter
groups  of $\Talg$,  respectively, both  isomorphic to  $\Z^r$ for  some $r$
called  the  rank  of $\Galg$  (or  of  $\Talg$  or  of $\Psi$).  These  are
in  duality  via  a  natural pairing  $\langle\cdot,\cdot\rangle:  X  \times
Y  \to  \Z$.  Here  $R$  is  a finite  subset  of  $X$,  called  the  roots.
There is  a bijection  ${}^\vee: R  \to R^\vee  \subset Y$,  $\alpha \mapsto
\alpha^\vee$, to  the set  $R^\vee$ of coroots,  such that  $\langle \alpha,
\alpha^\vee\rangle = 2$ for all $\alpha \in R$.

Each $\alpha  \in R$  defines reflections $s_\alpha:X\to  X$, $x\mapsto  x -
\langle  x,  \alpha^\vee  \rangle  \alpha$,  and  $s_\alpha^\vee:  Y\to  Y$,
$y\mapsto y-\langle\alpha, y\rangle \alpha^\vee$. The group $W$ generated by
all  $s_\alpha$  is called  the  Weyl  group of  $\Galg$  or  $\Psi$, it  is
isomorphic to the  group $W^\vee$ generated by the  $s_\alpha^\vee$. We have
$R W = R$ and $R^\vee W^\vee = R^\vee$.

Let $\Delta=\{\alpha_1,\ldots,\alpha_l\}  \subset R$   be  a  set  of simple
roots, that is each root is a  linear combination of simple ones with either
non-negative or  non-positive coefficients.  The integer  $l$ is  called the
semisimple rank of $\Galg$ and  $\Psi$. The set $\{s_\alpha\mid\; \alpha \in
\Delta\}$ is  a set of  Coxeter generators of  $W$ and $\Delta$  is linearly
independent as  subset of  $X\otimes_\Z \Q$. The  matrix $C  \in \Z^{l\times
l}$, $C_{ij}  = \langle\alpha_j,\alpha_i^\vee\rangle$ for $1\leq  i,j\leq l$
is called the  Cartan matrix of $\Psi$.  We have $\Delta W =  R$. The matrix
$C$ is the Cartan matrix of a crystallographic root system, the set $\Delta$
can be  reordered such  that $C$  has a block  diagonal form  whose diagonal
blocks are  in the  list given in~\cite[3.6]{Ca72}.  Cartan matrices  can be
encoded  in  a compact  way  by  Dynkin  diagrams,  this is  also  explained
in~\cite[3.6]{Ca72}.

We now introduce  a compact description of  a root datum which  is useful to
specify a root datum  and for computations. This is for  example used in the
GAP~\cite{GAP3} programs of the CHEVIE~\cite{CHEVIE} project.

Given $\Psi  = (X, R, Y,  R^\vee)$ we can  choose $\Z$-bases of $X$  and $Y$
which are dual to  each other and represent elements of $x\in  X$ and $y \in
Y$ by their coordinate row vectors with  respect to these bases (so, we have
$\langle x, y \rangle = y x^\tr$, where $^\tr$ means the transpose).

For $\Delta$  as above  we define  matrices $A,  A^\vee \in  \Z^{l\times r}$
where the $i$-th  row of $A$ contains the coordinates  of $\alpha_i$ and the
$i$-th row of $A^\vee$ those of $\alpha_i^\vee$.

From  $A$ and  $A^\vee$ we  can  compute the  whole root  datum: the  $i$-th
rows  of  the  two  matrices determine  the  generators  $s_{\alpha_i}$  and
$s_{\alpha_i}^\vee$ of $W$  and $W^\vee$, and the orbits of  the rows of $A$
under $W$  yield $R$ (and similarly  for $R^\vee$). The product  $C = A^\vee
A^\tr \in \Z^{l\times l}$ is the Cartan matrix of $\Psi$.

Vice  versa,  let $A,  A^\vee  \in  \Z^{l\times  r}$  be two  matrices  such
that  $C =  A^\vee A^\tr  \in  \Z^{l\times l}$  is  the Cartan  matrix of  a
crystallographic root  system, and let  $\bar k$ be an  algebraically closed
field. Then there exists a connected reductive algebraic group over $\bar k$
which  yields $(A,  A^\vee)$  as described  above  (use \cite[7.4.1,  9.5.1,
10.1]{Sp98}).

\begin{definition}\rm
We call  a pair of  matrices $(A,A^\vee) \in (\Z^{l\times  r})^2$ \emph{root
datum  matrices} if  $C =  A^\vee A^\tr  \in \Z^{l\times  l}$ is  the Cartan
matrix of a crystallographic root system.
\end{definition}

\begin{remark}\label{remscad}\rm
(a) Fixing the type  of a root datum via a Cartan  matrix $C \in \Z^{l\times
l}$  (or,  equivalently,  a  Dynkin diagram),  the  corresponding  connected
reductive groups  of adjoint type are  described by the root  datum matrices
$(\Id_l, C)$  (the simple roots are  a basis of $X$),  and the corresponding
groups  of simply-connected  type  are described  by  $(C^\tr, \Id_l)$  (the
simple coroots are a basis of $Y$).

(b) For $i = 1,2$ let $\Galg_i$ be a connected reductive group over $\bar k$
with a maximal  torus $\Talg_i$. Let $(A_i,A_i^\vee)$  be corresponding root
datum  matrices. Then  the direct  product $\Galg_1  \times \Galg_2$  can be
described with respect to the maximal torus $\Talg_1 \times \Talg_2$ by root
datum matrices $(A,A^\vee)$  where $A$ and $A^\vee$ are  block diagonal with
diagonal blocks $A_1, A_2$ and $A_1^\vee, A_2^\vee$, respectively.
\end{remark}

The following observation will be useful  later. We formulate it with roots,
there is a similar statement for the coroots.
\begin{lemma}\label{CequivLincomb}
Let  $\Psi = (X,R,Y,R^\vee)$ be a root datum and $\Delta \subset R$ a set of
simple roots.
\begin{itemize}
\item[(a)]  The   Cartan  matrix  $C  =   (\langle  \alpha_j,  \alpha_i^\vee
\rangle)_{i,j}$ or,  equivalently, the Dynkin  diagram of $\Psi$  labeled by
$\Delta$, determines the set $R$ of roots as linear combinations of those in
$\Delta$.
\item[(b)]  The  set  of  roots  $R$  as  linear  combinations  of  $\Delta$
determines the  Cartan matrix  $C$ or, equivalently,  the Dynkin  diagram of
$\Psi$.
\end{itemize}
\end{lemma}

\textbf{Proof.} (a) The set $R$ is the union of orbits of $\Delta$ under the
Weyl  group $W$  which  is  generated by  the  $s_\alpha$  with $\alpha  \in
\Delta$. For $\beta \in R$ (which  is a $\Z$-linear combination of $\Delta$)
we have $s_\alpha(\beta) =  \beta - \langle\beta, \alpha^\vee\rangle\alpha$,
so  the action  of $s_\alpha$  on the  $\Z$-lattice spanned  by $\Delta$  is
completely determined by the Cartan matrix.

(b) For any  two simple roots $\alpha_i, \alpha_j \in  \Delta$ the subset of
$R$ consisting  of linear combinations  of $\alpha_i$ and $\alpha_j$  is the
same as the  sub-root system spanned by these two  roots (see~\cite[Prop. in
1.10]{Hu90}).

So, to find the bond between $\alpha_i$ and $\alpha_j$ in the Dynkin diagram
labeled by $\Delta$ we  look at the subset of positive roots  in R which are
linear combinations of $\alpha_i$ and $\alpha_j$. There are $2$, $3$, $4$ or
$6$ such  roots, corresponding to  no, a single, a  double or a  triple bond
(types $2A_1$, $A_2$, $B_2$, $G_2$), respectively.  In the last two cases an
arrow must be added  pointing to the shorter root, this is the one occurring
with the largest coefficient in the linear combinations. \hfill$\Box$

\subsection{Homomorphisms of root data}
We recall  some information from  \cite[II 1.13--1.15]{Ja03}. For $i  = 1,2$
let  $\Galg_i$ be  connected reductive  groups over  the same  algebraically
closed field  $\bar k$, with  maximal tori $\Talg_i$ and  corresponding root
data $\Psi_i = (X_i, R_i, Y_i, R_i^\vee)$.

A homomorphism from  $\Psi_1$ to $\Psi_2$ is given by  a $\Z$-linear map $f:
X_2\to X_1$ such that $f$ induces a bijection $R_2 \to R_1$ and its dual map
$f^\vee: Y_1 \to Y_2$ induces a bijection $R_1^\vee \to R_2^\vee$.

For  each   such  homomorphism  of   root  data  there  is   a  homomorphism
$\phi:  \Galg_1 \to  \Galg_2$  that  maps $\Talg_1  \to  \Talg_2$ such  that
$\phi\mid_{\Talg_1}$  induces $f^\vee:  Y_1  \to Y_2$  and $\ker{\phi}  \leq
Z(\Galg_1) \leq  T_1$, where $Z(\Galg_1)$  denotes the center  of $\Galg_1$.
More precisely, $\ker{\phi} = \{t \in T_1\mid\; f(x)(t) = 1 $ for all $x \in
X_2\}$.

The map $\phi$ is surjective if and only if $Y_2/f^\vee(Y_1)$ is finite, and
$\phi$ is an isomorphism if and only if $f^\vee$ (or $f$) is invertible.

Moreover, $\phi$ is called  an isogeny if it is surjective  and has a finite
kernel, that is $f^\vee$ maps $Y_1$ injectively onto a finite index subgroup
of $Y_2$.

The root  datum associated to  a connected reductive  group is unique  up to
isomorphism.

To construct homomorphisms of root data we will use the following lemma.
\begin{lemma}\label{homlemma}
A homomorphism of root data is  determined by a $\Z$-linear map $f^\vee: Y_1
\to Y_2$  which induces a bijection  $R_1^\vee \to R_2^\vee$ and  which maps
$R_1^\perp$ to  $R_2^\perp$, where  $R_i^\perp =  \{y \in  Y_i\mid\; \langle
\alpha, y\rangle = 0$ for $\alpha \in R_i\}$.
\end{lemma}
\textbf{Proof.}
The given map $f^\vee: Y_1 \to Y_2$  induces its dual map $f:X_2 \to X_1$ as
follows. For $x_2\in  X_2$ the image $f(x_2) \in X_1$  is the unique element
such that $\langle f(x_2), y_1\rangle = \langle x_2, f^\vee(y_1)\rangle$ for
all $y_1 \in Y_1$.

We need to show that the map $f$ induces a bijection from $R_2 \to R_1$.

Let $\Delta_1 =\{\alpha_1,  \ldots, \alpha_l\}$ be a set of  simple roots in
$R_1$ and $\Delta_1^\vee =  \{\alpha_1^\vee, \ldots, \alpha_l^\vee\}$ be the
corresponding coroots. Since $f^\vee$ is $\Z$-linear and induces a bijection
$R_1^\vee \to  R_2^\vee$, it  must map  $\Delta_1^\vee$ to  a set  of simple
coroots of  $ R_2^\vee$.  So, there  is a  set of  simple roots  $\Delta_2 =
\{\beta_1, \ldots,  \beta_l\}$ of  $R_2$ such that  $f^\vee(\alpha_j^\vee) =
\beta_j^\vee$ for $1 \leq j \leq l$.

Now we use lemma~\ref{CequivLincomb}(b) to conclude that the Cartan matrices
of $\Psi_1$  and $\Psi_2$  are the same,  more precisely  $\langle \alpha_j,
\alpha_i^\vee\rangle  = \langle  \beta_j, \beta_i^\vee  \rangle$ for  all $1
\leq i, j \leq l$.

We show the lemma by checking that $f(\beta_i) = \alpha_i$ for $1\leq i \leq
l$.

Note that $\Q Y_1 = \Q \Delta_1^\vee \oplus \Q R_1^\perp$ because the Cartan
matrix of  $\Psi_1$ has full  rank $l$. So, we  can show that  $f(\beta_i) =
\alpha_i$  by  showing that  $\langle  f(\beta_i),  \alpha_j^\vee \rangle  =
\langle  \alpha_i, \alpha_j^\vee\rangle$  for $1  \leq  j \leq  l$ and  that
$\langle f(\beta_i), y\rangle  = \langle \alpha_i, y\rangle$ for  all $y \in
R_1^\perp$.

The first follows  because the Cartan matrices of $\Psi_1$  and $\Psi_2$ are
the  same:  $\langle  f(\beta_i), \alpha_j^\vee\rangle  =  \langle  \beta_i,
f^\vee(\alpha_j^\vee)\rangle  =  \langle  \beta_i,  \beta_j^\vee  \rangle  =
\langle \alpha_i, \alpha_j^\vee\rangle$.

The second follows because $f^\vee(y)  \in R_2^\perp$ for $y \in R_1^\perp$:
$\langle  f(\beta_i), y\rangle  = \langle  \beta_i, f^\vee(y)\rangle  = 0  =
\langle \alpha_i, y\rangle$. \hfill$\Box$

Of course, there is also a similar version of the lemma where the  roles  of
$X_i$ and $Y_i$ are interchanged.

\subsection{Frobenius morphisms}\label{ssec:frob}
From  now we  assume that  our field  $\bar k  = \bar\F_p$  is an  algebraic
closure of  the finite prime  field with $p$  elements, and that  $\Galg$ is
defined over  the finite subfield $\F_q  \leq \bar k$ with  $q$ elements. We
refer to~\cite[Chapter 3]{DM91} for an  explanation of this notion. There is
a corresponding  Frobenius morphism $F:\Galg\to\Galg$. We  consider the root
datum of $\Galg$  with respect to a maximal torus  $\Talg$ that is contained
in a Borel subgroup $\Balg$ with  $F(\Balg) = \Balg$ and $F(\Talg) = \Talg$.
Then $F$  induces a  map on $X$  which is  of the form  $q F_0$  where $F_0$
defines an automorphism of root data  of finite order which permutes the set
of simple  roots $\Delta$  that  is  determined   by  $\Balg$.  This follows
from~\cite[3.17]{DM91}  (the  $\tau$  in  that theorem  is  our  $F_0^{-1}$)
and~\cite[3.6(ii)]{DM91} (which shows that $F_0$ has finite order).

Vice versa,  each $q  F_0$ with  $F_0$ of  finite order  is induced  by some
Frobenius  morphism $F$  of $\Galg$  as  above; $F$  is uniquely  determined
by  $F_0$  and  $q$  up  to  conjugation  by  an  element  in  $\Talg$.  See
\cite[9.6]{Sp98} for more details.

The finite groups of fixed points  $G(q) = \Galg^F$ are called finite groups
of Lie type.  The group $G(q)$ is  determined up to isomorphism  by the root
datum $\Psi$ of $\Galg$, $F_0$ and $q$. (But various such tuples of data can
yield isomorphic groups $G(q)$.)

If the root  datum is described by root datum  matrices $(A,A^\vee)$ and the
elements of  $X$ and  $Y$ are considered  as row vectors  then $F_0$  can be
described by an invertible matrix in $\Z^{r\times r}$ of finite order.

We remark  that in this  setup we  do not cover  the Suzuki and  Ree groups.
These are fixed points of simple  reductive groups of types $B_2$, $F_4$ and
$G_2$  under generalized  Frobenius morphisms  whose square  is a  Frobenius
morphism as considered above (for $q$ an  odd power of $3$ in case $G_2$ and
an  odd  power  of  $2$  in  the  other  two  cases).  But  in  these  cases
parameterizations of the irreducible defining characteristic representations
are known, see~\ref{thm:finsc} and~\ref{SuzReeMain}.

\subsection{A covering group}\label{sec:covering}
A semisimple  group $\Galg$ has a  covering by a simply  connected group. In
this subsection we  explicitly construct such a covering  $\tilde \Galg$ for
general  connected reductive  $\Galg$. If  $F$  is a  Frobenius morphism  on
$\Galg$ we also construct a Frobenius  morphism $\tilde F$ on $\tilde \Galg$
which induces $F$ on $\Galg$.

\begin{proposition}\label{prop:covering}
Let  $\Galg$  be a  connected  reductive  group,  defined over  $\F_q$  with
Frobenius morphism $F$. Let the  root datum $\Psi=(X,R,Y,R^\vee)$ of $\Galg$
and $F$ be described by root datum matrices $(A,A^\vee)$ and $F_0$ as above.

There  are  root   datum  matrices  $(\tilde  A,  \tilde   A^\vee)$  and  an
automorphism $\tilde  F_0$ of finite  order of the corresponding  root datum
$\tilde  \Psi =  (\tilde  X, \tilde  R,  \tilde Y,  \tilde  R^\vee)$, and  a
homomorphism $\tilde \Psi \to \Psi$ with the following properties.
\begin{itemize}
\item[(a)] The  connected reductive  group $\tilde  \Galg$ over  $\bar \F_p$
determined by $\tilde  \Psi$ is a direct product of  simple simply connected
groups and a central torus $\tilde \Zalg^0$.
\item[(b)] The homomorphism $\tilde \Psi \to \Psi$ induces an isogeny $\pi: 
\tilde \Galg \to \Galg$.                                                    
\item[(c)]  $\pi$  induces  an  isomorphism from  $\tilde  \Zalg^0$  to  the
connected center $\Zalg^0$ of $\Galg$.
\item[(d)]  There   is  a  Frobenius  morphism   $\tilde  F:  \tilde\Galg\to
\tilde\Galg$  corresponding to  $\tilde F_0$  and $q$  which induces  $F$ on
$\Galg$.
\end{itemize}
\end{proposition}

\textbf{Proof.}  We first  construct  $\tilde A$  and  $\tilde A^\vee$.  Let
$C=A^\vee A^\tr \in  \Z^{l\times l}$ be the Cartan matrix  of $\Psi$ and $r$
be the rank of $\Psi$. Let $\tilde  A = (C^\tr\mid 0) \in \Z^{l\times r}$ be
the matrix  with $C$ as  the first $l$ columns  and $r-l$ zero  columns, and
similarly  let $\tilde  A^\vee =  (\Id_l\mid  0) \in  \Z^{l\times r}$.  Then
$(\tilde A, \tilde  A^\vee)$ are root datum matrices  because $\tilde A^\vee
\tilde A^\tr  = C$,  so determine  a root  datum $\tilde  \Psi =  (\tilde X,
\tilde R, \tilde Y, \tilde R^\vee)$  and a connected reductive group $\tilde
\Galg$ over $\bar\F_p$. After reordering of  the simple roots $\tilde A$ and
$\tilde A^\vee$  have block diagonal  form, the blocks corresponding  to the
simple components of  $\tilde \Galg$. So, it  is the root datum  of a direct
product of  the simple  components and  a torus.  The simple  components are
simply-connected, see remark~\ref{remscad}.

Let  $B   \in  \Z^{(r-l)\times   r}$  be  a   matrix  whose   rows  describe
a   $\Z$-basis   of  $R^\perp   \leq   Y$.   Then   the  matrix   $M^\tr   =
\left(\begin{array}{c}A^\vee\\B\end{array}\right)   \in    \Z^{r\times   r}$
describes  a  $\Z$-linear  map  $f^\vee:\tilde  Y \to  Y$  which  defines  a
homomorphism of root data: It maps  the simple coroots to simple coroots and
so by  lemma~\ref{CequivLincomb}(a) the coroots $\tilde  R^\vee$ to $R^\vee$
(the root data  have the same Cartan matrix). And  it induces an isomorphism
$\tilde R^\perp \to R^\perp$. Hence we can use Lemma~\ref{homlemma}.

We can  compute $B$  as follows:  Its rows are  a $\Z$-basis  of the  set of
solutions  $y \in  \Z^r$  of $y  A^\tr  =  0$. With  the  Smith normal  form
algorithm we can  compute invertible integer matrices $P$ and  $Q$ such that
$P A Q$ has diagonal form, so the  last $r-l$ columns of $A Q$ are zero (and
the first $l$  columns are $\Q$-linearly independent). We can  take the last
$r-l$ rows of $Q^\tr$ as matrix $B$.

The map  $f^\vee: \tilde Y  \to Y$ is injective,  its image is  generated by
$R^\vee$ and  $R^\perp$. So, the image  is invariant under $F_0^\tr$  and we
can  define $\tilde  F_0^\tr: \tilde  Y \to  \tilde Y$  by $\tilde  y \tilde
F_0^\tr := {f^\vee}^{-1} (f^\vee(\tilde y) F_0^\tr) = \tilde y M^\tr F_0^\tr
M^{-\tr}$. This  defines an automorphism  of finite order of  $\tilde \Psi$.
Now,  $\tilde\Psi$, $\tilde  F_0^\tr$  and  a prime  power  $q$ determine  a
reductive  $\tilde  \Galg$,  defined  over $\F_q$  with  Frobenius  morphism
$\tilde F$. We have a surjective homomorphism $\pi: \tilde \Galg \to \Galg$,
and $\tilde F$  induces a Frobenius morphism $F'$ on  $\Galg$, which induces
$F_0^\tr$  on  $Y$. So,  modifying  $\tilde  F$  by  a conjugation  with  an
appropriate  torus element  we can  assume that  $\tilde F$  induces $F$  on
$\Galg$.

The kernel $K := \ker(\pi)$ of the covering $\pi: \tilde \Galg \to \Galg$ is
finite because $M$ (and $M^\tr$) have  full $\Q$-rank $r$, so $XM\leq \tilde
X$ is of finite index. We have  $R^\perp = Y(\Zalg^0)$ and $\tilde R^\perp =
Y(\tilde\Zalg^0)$ and  since $f^\vee$  induces an isomorphism  between these
two lattices,  the homomorphism $\pi$ induces  an isomorphism $\tilde\Zalg^0
\to \Zalg^0$. In subsection~\ref{ssec:toruselts} we  show how to compute the
kernel of $\pi$ explicitly. 
\hfill$\Box$

\begin{lemma}\label{lemma:finfactor}
Let $\tilde\Galg$, $\Galg$  be connected reductive groups  with a surjective
homomorphism  $\pi:  \tilde\Galg \to  \Galg$  and  central kernel  $K$.  Let
$\tilde F$ be  a Frobenius morphism of $\tilde\Galg$ with  $\tilde F(K) = K$
and $F$ be the induced Frobenius  morphism on $\Galg$. The induced map $\pi:
\tilde\Galg^{\tilde F} \to \Galg^F$ is  in general not surjective. We define
$\cL(K) = \{z^{-1}\tilde F(z)\mid\; z\in K\}$. Then $\pi(\tilde\Galg^{\tilde
F})$ is a normal subgroup of $\Galg^F$ and there is a natural isomorphism
\[ \Galg^F/\pi(\tilde\Galg^{\tilde F}) \stackrel{\sim}{\to} K/\cL(K).\]
\end{lemma}
\textbf{Proof.}
We first  show that $\pi(\tilde\Galg^{\tilde  F})$ is normal. Let  $\tilde h
\in \tilde\Galg^{\tilde F}$, $g \in  \Galg^F$ and $\tilde g \in \tilde\Galg$
with $\pi(\tilde g) = g$. Then $\tilde F(\tilde g) = \tilde g z$ for some $z
\in K$. It  follows $\tilde F(\tilde g^{-1}) = z^{-1}  \tilde g^{-1}$. Hence
$\tilde F(\tilde g^{-1}  \tilde h \tilde g) = z^{-1}  \tilde g^{-1} \tilde h
\tilde g  z = \tilde  g^{-1} \tilde h  \tilde g$ because  $K$ and so  $z$ is
central. This shows that $g^{-1}\pi(\tilde h) g = \pi(\tilde g^{-1} \tilde h
\tilde g) \in \pi(\tilde\Galg^{\tilde F})$.

Since the group $K$ as subgroup of the center of $\tilde\Galg$ is abelian if
follows that $\cL(K)$ is a subgroup of $K$.

We have $\Galg \cong \tilde\Galg/K$ and  for $g\in \tilde \Galg$ we have $gK
\in (\tilde\Galg/K)^F \cong  \Galg^F$ if and only if  $g^{-1}\tilde F(g) \in
K$. We consider the map
\[  \Galg^F   \cong  (\tilde   \Galg/K)^F  \to  K/\cL(K),\quad   gK  \mapsto
g^{-1}\tilde F(g) \cL(K).\]
Since  $K$ is  central, it  is easy  to check  that this  is a  well-defined
homomorphism. The Lang-Steinberg theorem (for $\tilde\Galg$) shows that this
map is surjective. An element $gK$ is in  the kernel of this map if and only
if $gK$ contains an element of $\tilde \Galg^{\tilde F}$.
\hfill$\Box$

\subsection{Torus elements}\label{ssec:toruselts}
Given  a  root  datum  $\Psi  = (X,R,Y,R^\vee)$  for  $\Galg$  and  $\Talg$,
we  can  recover  $\Talg$  by   the  isomorphism  $\Talg  \cong  Y\otimes_\Z
\bar\F_p^\times$. Via some fixed  isomorphism we identify the multiplicative
group $\bar\F_p^\times$ with the additive  group $\Q_{p'}/\Z$ of elements of
$p'$-order in $\Q / \Z$. See \cite[3.1]{Ca85} for more details.

Choosing dual  bases of $X$  and $Y$, we can  describe $\Psi$ by  root datum
matrices $(A, A^\vee)$ and identify  $\Talg \cong Y \otimes_\Z (\Q_{p'}/\Z)$
with $r$-tuples  of elements in $\Q_{p'}/\Z$.  In this setup we  can compute
$y(c)$ for $y \in Y$ and $c \in  \Q_{p'}/\Z$, and apply $x \in X$ and $F$ to
$t \in \Talg = (\Q_{p'}/\Z)^r$ as follows:
\[ y(c) = c\cdot y, \quad\quad x(t) = t x^\tr \in \Q_{p'}/\Z, \quad\quad
F(t) = q t F_0^\tr \in \Talg.\]

The center  of $\Galg$ is  the intersection of  the kernels of  all (simple)
roots in $\Talg$.  We can compute it  as the solutions $t \in  \Talg$ of the
system of equations  $t A^\tr = 0 \in (\Q_{p'}/\Z)^l$.  The $F$-fixed points
$\Talg^F$  of $\Talg$  are the  solutions  $t \in  \Talg$ of  the system  of
equations $t (qF_0^\tr - \Id_r) = 0 \in (\Q_{p'}/\Z)^r$.

We consider  the isogeny from proposition~\ref{prop:covering},  $\pi: \tilde
\Galg  \to \Galg$.  In  the proof  of  the proposition  we  have computed  a
matrix $M$  describing the  map $f:  X \to \tilde  X$ for  the corresponding
homomorphism of root data.

We can compute the kernel $K$ of $\pi$  as set of solutions $t \in \Talg$ of
the system of equations
\[t M^\tr = 0 \in (\Q_{p'}/\Z)^r.\]
(The $\Z$-span of the rows of $M$ is the image $f(X) \leq \tilde X$.) And we
can compute the $\tilde F$-action on the elements $t\in K$ by
\[ \tilde F(t) = q t \tilde F_0^\tr.\]
This yields an  explicit description of the elements in  $K$, $K^{\tilde F}$
and $\cL(K)$.

\subsection{The derived subgroup}\label{ssec:Gprime}
We will also need to consider the  derived group $\Galg'$ of $\Galg$ and the
quotient torus $\Galg/\Galg'$. We  use the description in~\cite[8.1.9]{Sp98}
or~\cite[1.18]{Ja03}.

The images  of all coroots generate  $\Talg \cap \Galg'$ and  a character $x
\in  X$ has  $\Talg  \cap  \Galg'$ in  its  kernel if  and  only  if $x  \in
(R^\vee)^\perp$.

We compute a  matrix $D \in \GL_r(\Z)$  such that the last  $l-r$ columns of
$A^\vee  D$ are  zero  (for example  by the  Smith  normal form  algorithm).
Instead of  $(A, A^\vee)$  and $F_0$  we then  consider the  isomorphic data
$(AD^{-\tr}, A^\vee D)$ and $D^\tr F_0 D^{-\tr}$.

Now we get root datum matrices for  $\Galg'$ by taking the first $l$ columns
of $AD^{-\tr}$  and $A^\vee D$,  and the restriction  of $F$ to  $\Galg'$ is
described by the upper left $l\times l$ corner of $D^\tr F_0 D^{-\tr}$.

Furthermore,  the   lower  right   $(r-l)\times  (r-l)$  corner   of  $D^\tr
F_0  D^{-\tr}$  describes   the  Frobenius  action  induced   on  the  torus
$\Galg/\Galg'$.

\begin{lemma}\label{lem:quopiGprime}
Let $\pi: \tilde\Galg  = \tilde\Galg' \times Z^0 \to \Galg$  be the covering
and $\tilde F$ be the Frobenius  morphism of $\tilde\Galg$ as constructed in
proposition~\ref{prop:covering}.

Then $\pi(\tilde\Galg'^{\tilde F})  \leq \Galg^F$ is a  normal subgroup, and
the quotient $\Galg^F / \pi(\tilde\Galg'^{\tilde F})$ is an abelian group of
order prime to $p$.
\end{lemma}
\textbf{Proof.}
That $\pi(\tilde\Galg'^{\tilde F})  \leq \Galg^F$ is normal can  be shown as
in the  proof of  lemma~\ref{lemma:finfactor}, using that  $\tilde\Galg'$ is
normal in $\tilde\Galg$.

In~\cite[proof    of~13.20]{DM91}   it    is   shown    that   $\Galg^F    =
\Talg^F.\pi(\tilde\Galg'^{\tilde   F})$.  So,   the   quotient  $\Galg^F   /
\pi(\tilde\Galg'^{\tilde  F})$ is  isomorphic  to $\Talg^F  / (\Talg^F  \cap
\pi(\tilde\Galg'^{\tilde F}))$,  hence it is  abelian and of order  prime to
$p$.
\hfill$\Box$

\section{Irreducible representations in defining characteristic}
\label{sec:irrdefchar}

In  this  section  we  consider (finite  dimensional  rational)  irreducible
representations of our connected reductive  algebraic groups $\Galg$ and the
finite groups of Lie  type $\Galg^F$ over the defining field  $\bar k = \bar
\F_p$ of $\Galg$. As before, let $\Psi =(X,R,Y,R^\vee)$ be the root datum of
$\Galg$ and  $F_0:X\to X$ and $q$  be the finite order  automorphism and the
prime power determined by $F$.

\subsection{Representations of connected reductive groups}
In this subsection $\bar k$ can be  any algebraically closed field. We fix a
set $\Delta \subset R$ of simple roots. The set
\[ X_+ = \{x \in X\mid\; \langle x, \alpha^\vee \rangle \geq 0 \textrm{ for
} \alpha \in \Delta\} \]
is called the set of dominant weights of $\Galg$ (or $\Psi$).

One can associate  to each irreducible representation of  $\Galg$ over $\bar
k$ a highest weight $\lambda \in  X_+$. Chevalley proved the following basic
theorem, see~\cite[2.7]{Ja03}:

\begin{theorem}\label{thm:irralg}
Associating  the highest  weight induces  a bijection  from the  isomorphism
classes of irreducible  representations of $\Galg$ over $\bar k$  to the set
$X_+$ of dominant weights.
\end{theorem}

For  $\lambda  \in   X_+$  we  denote  by   $L(\lambda)$  the  corresponding
irreducible module, and by $\rho_\lambda$ the corresponding representation.

\subsection{Finite groups of Lie type, simply connected case}
In  this   subsection  we   assume  that  $\Galg$   is  semisimple   and  of
simply-connected type. In  this case the simple coroots are  a $\Z$-basis of
$Y$. The elements of the  dual basis $\{\omega_1, \ldots, \omega_l\} \subset
X$ are  called the fundamental weights.  For a positive integer  $b$ we call
the subset
\[ \begin{array}{rcl}
X_b & = &\{  x \in X_+\mid\; \langle x, \alpha^\vee\rangle  < b \textrm{ for
all } \alpha\in\Delta\} \\ & = & \{ a_1\omega_1 + \ldots + a_l\omega_l\mid\;
0 \leq a_i < b \textrm{ for } 1\leq i \leq l\}
\end{array}\]
of dominant weights the set of $b$-restricted weights.

Steinberg proved  the following theorem, see~\cite[13.3]{St68}.  Although we
have excluded the cases  of Suzuki and Ree groups from  our general setup we
include them in this theorem.

\begin{theorem}\label{thm:finsc}
(a) The restrictions  of $\rho_\lambda$ with $\lambda \in  X_q$ to $\Galg^F$
remain irreducible. This  induces a bijection from $X_q$  to the isomorphism
classes of irreducible representations of $\Galg^F$ over $\bar \F_q$.

(b) Let $\Galg$ be  of type $B_2$, $F_4$ or $G_2$, $q^2$ be  an odd power of
$p=2$, $2$ or $3$, respectively, and $F_0$  be of order $2$. We consider the
set $X'_q$ of dominant weights  $\sum_{i=1}^l a_i\omega_i$ with $0\leq a_i <
q\sqrt{p}$ if $\alpha_i$ is a short simple root and $0\leq a_i < q/\sqrt{p}$
otherwise. Then the  restrictions of $\rho_\lambda$ with  $\lambda \in X'_q$
induce a  bijection from  $X'_q$ to the  isomorphism classes  of irreducible
representations of $\Galg^F$ over $\bar \F_q$.
\end{theorem}

\subsection{Finite  tori}
Let  $\Galg  =  \Talg$  be  a  torus.  The  irreducible  representations  of
$\Talg$  are  the  characters  $X(T)$.  Let $F_0:  X\to  X$  be  the  finite
order  automorphism  induced   by  $F$,  and  let  $m  \in   \N$  such  that
$F_0^m  =  \operatorname{Id}$. We  have  $\Talg^F  \leq \Talg^{F^m}$.  Since
these  groups  are abelian,  each  irreducible  representation of  $\Talg^F$
can  be  extended  to  one  of  $\Talg^{F^m}$,  see~\cite[5.5]{Is76}.  Using
subsection~\ref{ssec:toruselts}   we  see   that  the   group  $\Talg^{F^m}$
is  isomorphic  to   a  direct  product  of  $r$  cyclic   groups  of  order
$q^m-1$.  And restriction  yields a  bijection  from the  set of  characters
$\{\rho_\lambda\mid\;  \lambda \in  X_{q^m-1}\}$ of  $\Talg$ to  the set  of
irreducible characters over  $\bar \F_q$ of the  finite group $\Talg^{F^m}$.

\begin{remark}\label{rem:fintori}\rm  
All irreducible representations of a  finite torus $\Talg^F$ over $\bar\F_q$
are restrictions  of irreducible  representations (characters) of  the torus
$\Talg$.
\end{remark}

In the general case there seems to be no nice description of a subset of $X$
which yields  the pairwise different  characters of $\Talg^F$, for  this one
has to compute  an explicit parameterization of $\Talg^F$. This  can be done
by  solving  the  system of  equations  $t  (q  F_0^\tr  - \Id_r)  =  0  \in
(\Q_{p'}/\Z)^r$,  as explained  in  subsection~\ref{ssec:toruselts} (we  see
that the order  of $\Talg^F$ is just the characteristic  polynomial of $F_0$
evaluated at $q$).

\subsection{Extending representations}\label{sec:extrep}

\begin{proposition}\label{prop:ext} 
Let $\Galg$ be a connected reductive  group over $\bar k$ ($=\bar\F_p$) with
Frobenius morphism $F$. Let $H \leq  \Galg^F$ be a normal subgroup such that
$\Galg^F/H$ is an abelian group of order prime to $p$. Then each irreducible
representation of $H$  over $\bar k$ can be extended  to a representation of
$\Galg^F$. And  each irreducible representation  of $\Galg^F$ over  $\bar k$
restricts irreducibly to $H$.
\end{proposition}
\textbf{Proof.} 
We use Clifford  theory, see for example~\cite[9.18]{Hu81}.  It follows that
the two  statements in the proposition  are equivalent. We show  the latter:
the restriction  of every irreducible  $\bar k\Galg^F$-module $V$ to  $H$ is
irreducible.

The restriction  is a direct sum of irreducible $\bar kH$-modules $W_i$, 
\[ V_H = \bigoplus_{i=1}^r W_i, \]
we show $r=1$.
Let $U$ be a Sylow-p-subgroup of $\Galg^F$.  Then $U \leq H$ because $H$ has
$p'$-index.  The  only  simple  $\bar  kU$-module  is  the  trivial  module.
Therefore, each $W_i$ must have at least a one-dimensional subspace on which
$U$ acts trivially. So, $V$ contains at least an $r$-dimensional subspace on
which $U$ acts trivially.

Now  we  use  that  the  group  $\Galg^F$  is  a  finite  group  with  split
$(B,N)$-pair in characteristic $p$, see~\cite[Cor. 4.2.5]{Ge03}. Thus we can
apply a result by Richen and Curtis that says that the subspace of $V$ fixed
by $U$ is one-dimensional, see~\cite[4.3(c)]{Cu70}. Hence $r=1$ and $V_H$ is
an irreducible $kH$-module.
\hfill$\Box$

\subsection{Parameterization of irreducible representations of finite
groups of Lie type}\label{sec:main}
We can now describe the main result of this paper.

As before, let $\Galg$ be a connected reductive group over $\bar k$, defined
over $\F_q$ with  corresponding Frobenius morphism $F$, given  by root datum
matrices  $(A,A^\vee)$ and  a finite  order  matrix $F_0$,  as explained  in
section~\ref{sec:rootdata}.

In  proposition~\ref{prop:covering} we  have  constructed  a covering  $\pi:
\tilde\Galg =  \tilde\Galg' \times Z^0  \to \Galg$ and a  Frobenius morphism
$\tilde F$  of $\tilde\Galg$ inducing $F$  on $\Galg$. We write  $K$ for the
kernel of $\pi$ and $\Galg'$ for the derived subgroup of $\Galg$.

\begin{theorem}\label{maintheorem}
The  irreducible  representations   of  $\Galg^F$  over  $\bar   k$  can  be
parameterized by the direct product of the following three sets
\begin{itemize}
\item[(A)]  the   $q$-restricted  weights   of  $\tilde\Galg'$   which  have
$K^{\tilde F} \cap \tilde\Galg'$ in their kernel,
\item[(B)] the group 
$%\Galg'^F / \pi(\tilde\Galg'^{\tilde F}) \cong 
K^{\tilde F} \cap \tilde\Galg'$, 
\item[(C)] and the group $(\Galg/\Galg')^F$.
\end{itemize}
\end{theorem}

\textbf{Proof.}
This    follows    from    Steinberg's    theorem~\ref{thm:finsc}    applied
to    $\tilde\Galg'^{\tilde   F}$    and    Clifford    theory,   see    for
example~\cite[9.18]{Hu81}. We give more details.

We know from lemma~\ref{lem:quopiGprime} that $\pi(\tilde\Galg'^{\tilde F})$
is  a normal  subgroup of  $\Galg^F$ with  abelian quotient  of order  prime
to  $p$.  Thus we  can  apply  proposition~\ref{prop:ext}  to see  that  all
irreducible $\bar  k\Galg^F$-modules are  extensions of irreducible  $\bar k
\pi(\tilde\Galg'^{\tilde F})$-modules. By Clifford  theory the extensions of
a  fixed  $\bar  k  \pi(\tilde\Galg'^{\tilde F})$-module  to  $\Galg^F$  are
parameterized  by the  group  of  linear characters  of  the quotient  group
$\Galg^F/\pi(\tilde\Galg'^{\tilde F})$  which is isomorphic to  the quotient
group itself.

The  irreducible representations  of $\pi(\tilde\Galg'^{\tilde  F})$ can  be
interpreted  as the  irreducible  representations $\tilde\Galg'^{\tilde  F}$
which  have $K^{\tilde  F} \cap  \tilde\Galg'$ in  their kernel.  And, since
$\tilde\Galg'$  is  simply  connected, the  irreducible  representations  of
$\tilde\Galg'^{\tilde  F}$  are   by  theorem~\ref{thm:finsc}  parameterized
by  the  $q$-restricted  weights  of   $\tilde\Galg'$.  An  element  $z  \in
Z(\tilde\Galg') \leq  \tilde T \cap \tilde\Galg'$  lies in the kernel  of an
irreducible representation with highest weight $\lambda\in \tilde X(\tilde T
\cap \tilde\Galg')$  if its  (only) eigenvalue is  $1$. This  eigenvalue can
be  read  off at  the  weight  space of  the  highest  weight by  evaluating
$\lambda$  at  $z$.  This  shows that  the  irreducible  representations  of
$\pi(\tilde\Galg'^{\tilde F})$ can be parameterized by the set~(A).

We can parameterize the  characters $\Hom(\Galg^F / \pi(\tilde\Galg'^{\tilde
F}),  \bar   k^\times)$  in   two  steps,  first   by  the   restriction  to
$\Hom(\Galg'^F   /  \pi(\tilde\Galg'^{\tilde   F}),   \bar  k^\times)$   and
then   by  the   characters  $\Hom(\Galg^F   /  \Galg'^F,   \bar  k^\times)$
(again  by   Clifford  theory  because  all   characters  in  $\Hom(\Galg'^F
/  \pi(\tilde\Galg'^{\tilde  F}),  \bar  k^\times)$  extend  to  $\Galg^F  /
\pi(\tilde\Galg'^{\tilde  F})$). The  latter  yields  our parameter  set~(C)
using  $\Galg^F/\Galg'^F  \cong  (\Galg/\Galg')^F$ which  follows  from  the
Lang-Steinberg theorem. The set~(B) we  get from the isomorphism $\Galg'^F /
\pi(\tilde\Galg'^{\tilde  F}) \cong  K^{\tilde F}  \cap \tilde\Galg'$.  This
follows from  lemma~\ref{lemma:finfactor} applied  to the induced  map $\pi:
\tilde\Galg'  \to  \Galg'$  which  has  the  finite  kernel  $K'  =  K  \cap
\tilde\Galg'$.  The lemma  shows  $\Galg'^F  / \pi(\tilde\Galg'^{\tilde  F})
\cong K'/\mathcal{L}(K')$. This last group is isomorphic to $K'^{\tilde F} =
K^{\tilde  F} \cap  \tilde\Galg'$  which follows  from  dualizing the  exact
sequence $1 \to K'^{\tilde F} \to K' \to \mathcal{L}(K') \to 1$.
\hfill$\Box$

We  now indicate  how to  compute the  parameter sets~(A),  (B) and~(C).  In
proposition~\ref{prop:covering}  we  have  constructed  the  root  datum  of
$\tilde\Galg$  such  that the  first  $l$  coordinates  and the  last  $r-l$
coordinates of  $\tilde X$ and $\tilde  Y$ correspond to the  factors of the
direct  product  $\tilde T  =  (\tilde  T  \cap \tilde\Galg')  \times  Z^0$.
Thus  it is  easy  to  decide which  elements  of  $K^{\tilde F}$,  computed
as  in  subsection~\ref{ssec:toruselts},  are contained  in  $\tilde\Galg'$,
this yields the set~(B).

In~\ref{ssec:toruselts} we  have also shown  how to evaluate a  $\lambda \in
X(\tilde T  \cap \tilde\Galg')$ at a  torus element. This way  we can decide
which  $q$-restricted  weights of  $\tilde\Galg'$  have  $K^{\tilde F}  \cap
\tilde\Galg'$ in their kernel. This determines the set~(A).

For  the set~(C)  we need  to  compute the  structure of  the abelian  group
$\Galg^F/\Galg'^F  \cong (\Galg/\Galg')^F$.  In subsection~\ref{ssec:Gprime}
we have described how to compute the $F$-action on the torus $\Galg/\Galg'$.
We can  use~\ref{ssec:toruselts} again  to compute  the $F$-fixed  points of
$\Galg/\Galg'$.

\begin{remark}\label{rem:restralg}\rm
(a) Assume that  the derived group $\Galg'$ of $\Galg$  is simply connected.
Then  each irreducible  $\bar k  \Galg^F$-module  is the  restriction of  an
irreducible  $\bar k  \Galg$-module. In~\cite[App.  1.3]{Hz09} Herzig  gives
another parameterization in this case:  Namely by all $q$-restricted weights
of $\Galg$  (these are  infinitely many  if $\Galg$  is not  semisimple) and
showing  that two  $q$-restricted weights  $\lambda_1, \lambda_2$  yield the
same  restriction  to $\Galg^F$  if  and  only if  $\lambda_1-\lambda_2  \in
(q\cdot \textrm{id} - F_0) (R^\vee)^\perp$.

(b)  In   general,  not  all   irreducible  $\bar  k   \Galg^F$-modules  are
restrictions  of modules  of  the  algebraic group  $\Galg$.  As an  example
consider  $\Galg=\PGL_{l+1}(\bar{k})$, the  adjoint  groups  of type  $A_l$,
with  Frobenius  map  $F$  such   that  $\Galg^F  =  \PGL_{l+1}(q)$.     For
some  prime   powers  $q$  the   finite  group  $\Galg^F$   has  non-trivial
$\bar{k}$-representations  of dimension  $1$. Such  representations are  not
restrictions from $\Galg$ because $\Galg$ is perfect.
\end{remark}
\textbf{Proof.}  We show  the first  statement  of~(a) using  our setup.  If
$\Galg'$  is  simply connected  our  parameterization  of irreducible  $\bar
k\Galg^F$ modules  in theorem~\ref{maintheorem} is particularly  simple: the
set~(A) consists of all $q$-restricted weights of $\Galg'$, the group~(B) is
trivial and~(C) is the finite torus $(\Galg/\Galg')^F$.

Since  $\Galg'$ is  simply connected,  $X$ contains  $\tilde \omega_i$  with
$\langle \tilde\omega_i, \alpha_j^\vee \rangle = \delta_{ij}$ for $1\leq i,j
\leq l$. So, for each $q$-restricted  weight $\lambda'$ of $\Galg'$ there is
a $\lambda \in X$ such that  the module $L(\lambda)$ of $\Galg$ restricts to
$\Galg'$ as $L(\lambda')$. Together with Steinberg's theorem~\ref{thm:finsc}
this shows  that each irreducible representation  $\tilde\rho$ of $\Galg'^F$
can  be  extended  to  a   representation  $\rho$  of  the  algebraic  group
$\Galg$.  All  the  other  extensions   of  $\tilde\rho$  to  $\Galg^F$  are
obtained  by  tensoring  $\rho|_{\Galg^F}$  with the  linear  characters  of
$\Galg^F/\Galg'^F$. But  these are also  obtained as restrictions  of linear
characters  of  the  algebraic  group  $\Galg/\Galg'$ as  we  have  seen  in
remark~\ref{rem:fintori}.
\hfill$\Box$

\subsection{A variant}\label{ssec:variant}
As  a variant  of  theorem~\ref{maintheorem}  we could  have  first given  a
parameterization of the  irreducible representations of $\tilde\Galg^{\tilde
F}$ and  use Clifford  theory only  for the  quotient $\Galg^F  / \pi(\tilde
\Galg^{\tilde F})$.  But our description in  theorem~\ref{maintheorem} often
leads to a more natural parameterization.

For example,  let $\Galg = \GL_{l+1}(\bar  k)$ and $q \equiv  1 \pmod{l+1}$.
Then $\tilde\Galg =  \SL_{l+1}(\bar k) \times Z^0$ and the  kernel of $\pi$,
$K  =  K^{\tilde  F}$,  is  cyclic  of order  $l+1$  and  is  isomorphic  to
$\Galg^F / \pi(\tilde\Galg^{\tilde F})$.  The irreducible representations of
$\tilde\Galg^{\tilde  F} \cong  \SL_{l+1}(q) \times  (Z^0)^{\tilde F}$  (the
second factor is  cyclic of order $q-1$)  are easy to describe. But  it is a
bit complicated  to describe  the subset  which has $K$  in its  kernel. The
quotient $\Galg^F / \pi(\tilde \Galg^{\tilde  F})$ is cyclic of order $l+1$,
so its irreducible representations are also easy to describe.

Our  parameterization  in  theorem~\ref{maintheorem}   is  more  natural  in
this  example: The  derived  subgroup  of $\Galg$  and  of $\tilde\Galg$  is
$\SL_{l+1}(\bar  k)$  and so  is  simply  connected.  Hence  we are  in  the
situation  of  remark~\ref{rem:restralg}(a),  our set~(A)  consists  of  all
$q$-restricted  dominant  weights of  $\SL_{l+1}(\bar  k)$  and our  set~(B)
is  trivial.  The  set~(C)  corresponds to the $q-1$  linear  characters  of
$\GL_{l+1}(q)/\SL_{l+1}(q)$.

\begin{remark}\label{SuzReeMain}\rm
A variant  of the main  theorem~\ref{maintheorem} is also true  if $\Galg^F$
has Suzuki or Ree groups as components. Since the Suzuki and Ree groups have
trivial center, we can assume that  $\Galg^F$ arises from an algebraic group
such that  the Suzuki and Ree  components are coming from  direct factors of
$\Galg$ of  simply connected type.  We can  then deal with  these components
using theorem~\ref{thm:finsc}(b).
\end{remark}

\subsection{An example}\label{ssec:example}
Let us consider as an example a  reductive group $\Galg$ which occurs as the
centralizer  of  a semisimple  element  in  the  simple algebraic  group  of
type  $E_8$,  equipped  with  a  Frobenius morphism  $F$.  It  is  given  by
root  datum matrices  $A$, $A^\vee$,  and a  matrix $F_0$,  as explained  in
sections~\ref{ssec:rootdata} and~\ref{ssec:frob}:
\[
A := 
{\tiny
\left(\begin{array}{rrrrrrrr}%
1&0&0&0&0&0&0&0\\%
0&1&0&0&0&0&0&0\\%
0&0&1&0&0&0&0&0\\%
0&0&0&0&1&0&0&0\\%
0&0&0&0&0&1&0&0\\%
0&0&0&0&0&0&0&1\\%
2&3&4&6&5&4&3&1\\%
\end{array}\right)},\quad
\!\! A^\vee\!\! :=
{\tiny
\left(\begin{array}{rrrrrrrr}%
2&0&-1&0&0&0&0&0\\%
0&2&0&-1&0&0&0&0\\%
-1&0&2&-1&0&0&0&0\\%
0&0&0&-1&2&-1&0&0\\%
0&0&0&0&-1&2&-1&0\\%
0&0&0&0&0&0&-1&2\\%
0&0&0&0&0&0&1&-1\\%
\end{array}\right)},
\]
\[
F_0 :=
{\tiny 
\left(\begin{array}{rrrrrrrr}%
0&0&1&0&0&0&0&0\\%
0&1&0&0&0&0&0&0\\%
1&0&0&0&0&0&0&0\\%
-1&-1&-1&-1&-1&-1&-1&-1\\%
0&0&0&0&0&0&0&1\\%
2&3&4&6&5&4&3&1\\%
-2&-3&-4&-6&-5&-3&-2&-1\\%
0&0&0&0&1&0&0&0\\%
\end{array}\right)}.
\]

We do not fix  the $q$, but we want to investigate all  finite groups of Lie
type for any prime power $q$ which are determined by these data.

The group $G$  has rank $8$ and  semisimple rank $7$. Looking  at the Cartan
matrix $A^\vee A^\tr$ we see that the  pairs of simple roots number 1 and 3,
number 4 and 5 and number 6 and 7 each span a sub-root system of type $A_2$,
and root number 2 spans a subsystem of type $A_1$. The matrix $A F_0$ yields
a  permutation of  the rows  of $A$,  the permutation  is $(1,3)(4,6)(5,7)$.
Thus,  the data  describe groups  $\Galg^F$  which are  central products  of
components of  type $^2A_2(q)$, $A_2(q^2)$,  $A_1(q)$ and a finite  torus of
rank $1$.

Now we  look at the covering  group $\tilde\Galg$ of $\Galg$  constructed in
proposition~\ref{prop:covering}. We do not need  the matrices $\tilde A$ and
$\tilde A^\vee$,  but the  matrix $M^\tr$ is  essential which  describes the
homomorphism $\tilde Y \to Y$ that determines the covering $\pi: \tilde\Galg
\to \Galg$. As described in  the proof of proposition~\ref{prop:covering} we
can compute  a $\Z$-basis  of $R^\perp$  by applying  the Smith  normal form
algorithm to $A$.  This yields invertible integer matrices $P$  and $Q$ such
that $PAQ$ is of  diagonal form (the diagonal entries are  six times $1$ and
one $3$). Then $M^\tr$ is given by the  rows of $A^\vee$ and the last row of
$Q^\tr$, the latter spans $R^\perp$.  We furthermore need $\tilde F_0$ which
defines the Frobenius  morphism on the covering group  $\tilde\Galg$. We can
compute it with $M^\tr$ as $\tilde F_0^\tr = M^\tr F_0^\tr M^{-\tr}$. We get
\[
M^\tr =
{\tiny
\left(\begin{array}{rrrrrrrr}%
2&0&-1&0&0&0&0&0\\%
0&2&0&-1&0&0&0&0\\%
-1&0&2&-1&0&0&0&0\\%
0&0&0&-1&2&-1&0&0\\%
0&0&0&0&-1&2&-1&0\\%
0&0&0&0&0&0&-1&2\\%
0&0&0&0&0&0&1&-1\\%
0&0&0&1&0&0&-2&0\\%
\end{array}\right)},\;
\tilde F_0 =
{\tiny
\left(\begin{array}{rrrrrrrr}%
0&0&1&0&0&0&0&0\\%
0&1&0&0&0&0&0&0\\%
1&0&0&0&0&0&0&0\\%
0&0&0&0&0&1&0&0\\%
0&0&0&0&0&0&1&0\\%
0&0&0&1&0&0&0&0\\%
0&0&0&0&1&0&0&0\\%
0&0&0&0&0&0&0&1\\%
\end{array}\right)}.
\]

Using  these two  matrices we  can determine  the finite  kernel $K$  of the
covering  $\pi:  \tilde\Galg \to  \Galg$  and  its $\tilde  F$-fixed  points
$K^{\tilde F}$, as explained in section~\ref{ssec:toruselts}. To find $K$ we
solve the system of equations
\[ tM^\tr = 0 \in (\Q_{p'}/\Z)^r. \]
To do so, we use again the  Smith normal form algorithm to find matrices $P,
Q \in \GL_r(\Z)$  such that $PM^\tr Q$ is diagonal.  The diagonal entries in
our example  are six times $1$,  $3$ and $6$. It  is easy to write  down the
solutions of
\[ t_1 (PM^\tr Q) = 0 \in (\Q_{p'}/\Z)^r. \]
If the $i$-th diagonal  entry of the diagonal matrix is  an integer $n$ then
the  $i$-th entry  of any  solution $t_1$  has the  form $i/n_{p'}$  for one
$0  \leq  i  <  n_{p'}$  (where  $n_{p'}$ is  the  largest  divisor  of  $n$
prime  to $p$).  Having  found all  solutions $t_1$  of  this last  equation
we  get  the  solutions  of  the  original equation  as  $t  =  t_1  P$.  In
practice  we first  compute  all  solutions $t  \in  (\Q/\Z)^r$, because  we
have  not  yet  said  anything  about  the  $q$  and  so  the  $p$.  In  our
example  we have  $18$  solutions  for $t_1$  over  $\Q/\Z$,  they have  the
form  $(0,0,0,0,0,0,\frac{i}{3},\frac{j}{6})$. And  multiplying with  $P$ we
get  for $t$  the  $18$  $\Z$-linear combinations  of  the  two elements  $(
\frac{2}{3},  0, \frac{1}{3},  \frac{1}{3}, \frac{2}{3},  0, 0,  \frac{2}{3}
)$  and $(\frac{1}{3},  \frac{1}{2}, \frac{2}{3},  \frac{1}{3}, \frac{2}{3},
\frac{2}{3}, \frac{1}{3}, \frac{1}{2})$.

We find $K^{\tilde F}$ by applying  the Frobenius $\tilde F$ to the elements
just found. This action is for any $t \in \tilde\Talg$ given by
\[\tilde F(t) = t (q \tilde F_0^\tr). \]
To be able to evaluate this on $t \in  K$ we need to know the residue of $q$
modulo  all  denominators  of  the  coordinates  of  $t  \in  K$.  A  common
denominator of all these entries is
\[ \begin{array}{rcl} 
m & := & \textrm{ the largest elementary divisor of } M^\tr \\
& = & \lcm(\textrm{entries of Smith normal form of } M^\tr) = 6.
 \end{array}
\]

We still  do not fix  $q$, but the remaining  computations are done  for any
congruence class $c$  of a prime power modulo $m$  separately, assuming that
$q \equiv c \pmod m$. In our example  we have to distinguish the cases of $q
\equiv 1,2,3,4,5 \pmod 6$.

In cases $c = 2$ or $4$ we  have $p=2$ (the prime dividing $c$ and $m$), and
in this case the kernel $K$ only contains the $9$ elements given above which
are of order  $1$ or $3$. Similarly, in  case $c = 3$ we have  $p=3$ and $K$
only contains the two elements of order  $1$ and $2$. In the other cases $K$
contains all $18$ elements given above.

For  the  computation of  $K^{\tilde  F}$  we  comment  on the  case  $c=2$.
Multiplying the  elements of  $K$ by  $q \tilde F_0^\tr$  and using  that $q
\equiv 2 \pmod 6$  we find that only the three  multiples of $( \frac{1}{3},
0, \frac{2}{3}, \frac{1}{3}, \frac{2}{3},  \frac{2}{3}, \frac{1}{3}, 0)$ are
mapped to themselves.

We  need  to  decide  which  of  these  $\tilde  F$-fixed  elements  lie  in
$\tilde\Galg'$. This is easy to see, because the first $l$ basis elements of
$\tilde X$ and $\tilde Y$ correspond  to the maximal torus of the semisimple
factor and  derived subgroup $\tilde\Galg'$.  So, here the  $\tilde F$-fixed
elements of  $K$ all lie in  $K^{\tilde F} \cap \tilde\Galg'$  because their
last coordinate is $0$.

Considering also  the other cases  for $c$ we  find that $K^{\tilde  F} \cap
\tilde\Galg'$ is cyclic of  order $3$ if $q \equiv 2$ or $5  \pmod 6$ and it
is  trivial  in  the other  cases.  We  have  found  the group~(B)  for  our
parameterization of the irreducible representations of $\Galg^F$.

The  elements  of  $K^{\tilde F} \cap \tilde\Galg'$ are  also needed to find
our parameter set~(A). This consists  of all $q$-restricted weights $\lambda
\in \tilde  X^+$ of $\tilde\Galg'$ which  are trivial on $K^{\tilde  F} \cap
\tilde\Galg'$. This means
\[t\lambda^\tr = 0 \in \Q_{p'}/\Z \textrm{ for all } t \in K^{\tilde F} \cap
\tilde\Galg'.\]
These equations can be reformulated in terms of integers by multiplying with
a  common multiple  $m'$ of  all denominators  in $t  \in K^{\tilde  F} \cap
\tilde\Galg'$ (a divisor of $m$). In our  example we can multiply with $m' =
3$ and then consider the equations modulo $m'$:
\[ (m' t) \lambda^\tr \equiv 0 \pmod {m'} \textrm{ for all } t \in K^{\tilde
F} \cap \tilde\Galg'.\]
Writing all  $(m' t)$  for a  set of  generators $t$  of $K^{\tilde  F} \cap
\tilde\Galg'$ in one matrix we can  further simplify the system of equations
by computing  the Hermite normal  form $\pmod {m'}$  of this matrix.  In our
example we get no non-trivial equation  if $c \not\in \{2,5\}$. So, in these
cases  all  $q^l$ $q$-restricted  weights  $\lambda$  lie in  our  parameter
set~(A). If $c  = 2$ or $5$, the set~(A)  contains only those $q$-restricted
$\lambda$ which fulfill the single equation
\[  ( 1, 0, 2, 1, 2, 2, 1 ) \lambda^\tr = 0 \pmod 3. \]

Using  this  equation   it  is  easy  to  check  for   a  concrete  $q$  and
$q$-restricted weight if it is in the parameter set~(A).

For general $q$ we  can also count the number of  parameters in the set~(A).
For this we use the following trivial lemma.

\begin{lemma}\label{lemma:count}
Let $q,c,i,m \in  \N$ with $0\leq i, c  < m$ and $q \equiv c  \pmod m$. Then
the number of integers $j$ with $0\leq j  < q-1$ and $j \equiv i \pmod m$ is
$(q-c)/m$ if $i \geq c$ and $(q-c)/m + 1$ for $i < c$.
\end{lemma}

This lemma can be applied recursively  to count the sets~(A). For example in
the  case $q  \equiv  2 \pmod  3$  above we  need to  count  the $\lambda  =
(\lambda_1,\ldots,  \lambda_7) \in  \Z^7$ with  $0 \leq  \lambda_i <  q$ for
$i=1,\ldots 7$ and $1 \cdot \lambda_1 + 0 \cdot \lambda_2 + \cdots + 1 \cdot
\lambda_7 \equiv 0 \pmod 3$.

From the lemma  we can easily deduce how often  each congruence class $\pmod
3$ is hit  by $1\cdot \lambda_1$, $0 \cdot \lambda_2$,  and so on. Combining
this it is easy to count how often each congruence class $\pmod 3$ is hit by
$1\cdot \lambda_1 + 0\cdot \lambda_2$. In  the next step we find the numbers
for  the  expressions  $1\cdot  \lambda_1  +  0\cdot  \lambda_2  +  2  \cdot
\lambda_3$. Going on  recursively, we find for each  congruence class $\pmod
3$ the  number of $q$-restricted  $\lambda$ with $(  1, 0, 2,  1, 2, 2,  1 )
\lambda^\tr$  in  that class.  In  particular,  we  find for  the  $0$-class
the  number  of  $q$-restricted  weights  in~(A),  it  is  $(q^7+2q)/3$  for
$q=2,5\pmod{6}$.

Finally, we  need the  set~(C), the  structure of  $(\Galg/\Galg')^F$. Using
subsection~\ref{ssec:Gprime} we can  find the matrix of $F_0$  acting on the
characters of this  torus via the transformation of the  matrix $A$ to Smith
normal form. In our example we find  the $1\times 1$ identity matrix. So the
group of $F$-fixed points in this torus is cyclic of order $q-1$. In general
the  order of  a  finite torus  is the  characteristic  polynomial of  $F_0$
evaluated at  $q$. The precise structure  of the finite abelian  group for a
specific $q$ is found by the  Smith normal form of the characteristic matrix
at $q$. See~\cite[Chapter 3]{Ca85} for more details.

To summarize: The parameter group~(C) is  for any $q$ cyclic of order $q-1$.
For $q  \equiv 2\pmod{3}$ the  parameter group~(B) is  of order $3$  and the
set~(A)  contains  $(q^7+2q)/3$  weights.  For $q  \equiv  0,1\pmod{3}$  the
group~(B)  is trivial  and  the set~(A)  contains  all $q^7$  $q$-restricted
weights.

\section{The case when $\Galg$ is simple}\label{sec:simple}
In  this   last  section   of  the   paper  we  want   to  apply   our  main
theorem~\ref{maintheorem}  to all  finite groups  of Lie  type arising  from
simple algebraic groups  $\Galg$. As an application we  determine the number
of semisimple classes of these groups.

As  before, we  exclude  here the  Suzuki  and Ree  groups,  in these  cases
the  $q^2$, respectively  $q^4$,  irreducible  representations were  already
described in theorem~\ref{thm:finsc}(b).

For each  type of  irreducible root system  $R$, we choose  a set  of simple
roots $\Delta=\{\alpha_1,\ldots,\alpha_l\}\subseteq  R$. We fix  a numbering
of the simple roots via the Dynkin diagrams given in table~\ref{tab:dynkin}.
The  node labelled  by  $i$ corresponds  to the  simple  root $\alpha_i$  of
$\Delta$.   This is the labelling used in \textsc{Chevie}; see~\cite{CHEVIE}
(the  often used  Bourbaki  labelling  is different  for  types  $B, C,  D$,
where it starts to count from the right side of the shown diagrams).

For a Frobenius morphism $F$ of $\Galg$  we consider a root datum of $\Galg$
with respect to a maximally split maximal torus. Then $F_0$ permutes the set
of simple roots and induces a graph automorphism of the Dynkin diagram. This
graph automorphism can  be non-trivial in cases $A_l$ with  $l\geq 2$, $D_l$
with $l\geq 4$ and $E_6$. We also write $F_\epsilon$ instead of $F$ in these
cases  with  $\epsilon  =1$  in  case of  the  trivial  graph  automorphism,
$\epsilon =  -1$ in case of  the graph automorphism of  order $2$ (permuting
nodes $1$  and $2$ in case  $D_l$) and $\epsilon=3$ in  case $D_4$ permuting
the nodes with cycle $(1,2,4)$.

Let $\Galg_{sc}$  be the simply-connected simple  group of the same  type as
$\Galg$. As explained in  proposition~\ref{prop:covering} we have an isogeny
$\Galg_{sc} \to  \Galg$ with  a central  kernel $K$  and $\Galg_{sc}$  has a
Frobenius morphism that  induces $F$ on $\Galg$, we denote  that also by $F$
or $F_\epsilon$.

As in the  proof of proposition~\ref{prop:covering} we choose  as root datum
matrices for  $\Galg_{sc}$ the pair  $(C^\tr,\Id)$, where $C$ is  the Cartan
matrix corresponding to the chosen numbering of $\Delta$. This means that in
the  root  datum  $(\tilde  X,  \tilde  R,  \tilde  Y,  \tilde  R^\vee)$  of
$\Galg_{sc}$  we use  the simple  coroots  as basis  of $\tilde  Y$ and  the
fundamental weights  as basis  of $\tilde  X$. The matrix  for $F_0$  is the
permutation matrix  for the  graph automorphism induced  by $F$.  As before,
we  identify a  maximal  torus  of $\Galg_{sc}$  with  $\tilde Y  \otimes_\Z
(\Q_{p'}/\Z) \cong (\Q_{p'}/\Z)^l$.

\begin{table}
\begin{center}
{
\newlength{\tmplen}
\setlength{\tmplen}{\unitlength}
\setlength{\unitlength}{0.8pt}
\scriptsize
\begin{picture}(360,190)
\put( 10, 40){$E_7$}
\put( 40, 40){\circle*{5}}
\put( 38, 45){1}
\put( 40, 40){\line(1,0){20}}
\put( 60, 40){\circle*{5}}
\put( 58, 45){3}
\put( 60, 40){\line(1,0){20}}
\put( 80, 40){\circle*{5}}
\put( 78, 45){4}
\put( 80, 40){\line(0,-1){20}}
\put( 80, 20){\circle*{5}}
\put( 85, 18){2}
\put( 80, 40){\line(1,0){20}}
\put(100, 40){\circle*{5}}
\put( 98, 45){5}
\put(100, 40){\line(1,0){20}}
\put(120, 40){\circle*{5}}
\put(118, 45){6}
\put(120, 40){\line(1,0){20}}
\put(140, 40){\circle*{5}}
\put(138, 45){7}

\put(190, 40){$E_8$}
\put(220, 40){\circle*{5}}
\put(218, 45){1}
\put(220, 40){\line(1,0){20}}
\put(240, 40){\circle*{5}}
\put(238, 45){3}
\put(240, 40){\line(1,0){20}}
\put(260, 40){\circle*{5}}
\put(258, 45){4}
\put(260, 40){\line(0,-1){20}}
\put(260, 20){\circle*{5}}
\put(265, 18){2}
\put(260, 40){\line(1,0){20}}
\put(280, 40){\circle*{5}}
\put(278, 45){5}
\put(280, 40){\line(1,0){20}}
\put(300, 40){\circle*{5}}
\put(298, 45){6}
\put(300, 40){\line(1,0){20}}
\put(320, 40){\circle*{5}}
\put(318, 45){7}
\put(320, 40){\line(1,0){20}}
\put(340, 40){\circle*{5}}
\put(338, 45){8}

\put( 10, 80){$G_2$}
\put( 40, 80){\circle*{5}}
\put( 38, 85){1}
\put( 40, 78){\line(1,0){20}}
\put( 40, 80){\line(1,0){20}}
\put( 40, 82){\line(1,0){20}}
\put( 46, 78){$>$}
\put( 60, 80){\circle*{5}}
\put( 58, 85){2}

\put(100, 80){$F_4$}
\put(130, 80){\circle*{5}}
\put(128, 85){1}
\put(130, 80){\line(1,0){20}}
\put(150, 80){\circle*{5}}
\put(148, 85){2}
\put(150, 78){\line(1,0){20}}
\put(150, 81){\line(1,0){20}}
\put(156, 77){$>$}
\put(170, 80){\circle*{5}}
\put(168, 85){3}
\put(170, 80){\line(1,0){20}}
\put(190, 80){\circle*{5}}
\put(188, 85){4}

\put(230, 80){$E_6$}
\put(260, 80){\circle*{5}}
\put(258, 85){1}
\put(260, 80){\line(1,0){20}}
\put(280, 80){\circle*{5}}
\put(278, 85){3}
\put(280, 80){\line(1,0){20}}
\put(300, 80){\circle*{5}}
\put(298, 85){4}
\put(300, 80){\line(0,-1){20}}
\put(300, 60){\circle*{5}}
\put(305, 58){2}
\put(300, 80){\line(1,0){20}}
\put(320, 80){\circle*{5}}
\put(318, 85){5}
\put(320, 80){\line(1,0){20}}
\put(340, 80){\circle*{5}}
\put(338, 85){6}

\put( 10,130){$D_l$}
\put( 40,150){\circle*{5}}
\put( 45,150){1}
\put( 40,150){\line(1,-1){21}}
\put( 40,110){\circle*{5}}
\put( 45,108){2}
\put( 40,110){\line(1,1){21}}
\put( 60,130){\circle*{5}}
\put( 58,135){3}
\put( 60,130){\line(1,0){30}}
\put( 80,130){\circle*{5}}
\put( 78,135){4}
\put(100,130){\circle*{1}}
\put(110,130){\circle*{1}}
\put(120,130){\circle*{1}}
\put(130,130){\line(1,0){10}}
\put(140,130){\circle*{5}}
\put(138,135){$l$}

\put(210,130){$C_l$}
\put(240,130){\circle*{5}}
\put(238,135){1}
\put(240,128){\line(1,0){20}}
\put(240,131){\line(1,0){20}}
\put(246,127){$>$}
\put(260,130){\circle*{5}}
\put(258,135){2}
\put(260,130){\line(1,0){30}}
\put(280,130){\circle*{5}}
\put(278,135){3}
\put(300,130){\circle*{1}}
\put(310,130){\circle*{1}}
\put(320,130){\circle*{1}}
\put(330,130){\line(1,0){10}}
\put(340,130){\circle*{5}}
\put(338,135){$l$}

\put( 10,170){$A_l$}
\put( 40,170){\circle*{5}}
\put( 38,175){1}
\put( 40,170){\line(1,0){20}}
\put( 60,170){\circle*{5}}
\put( 58,175){2}
\put( 60,170){\line(1,0){30}}
\put( 80,170){\circle*{5}}
\put( 78,175){3}
\put(100,170){\circle*{1}}
\put(110,170){\circle*{1}}
\put(120,170){\circle*{1}}
\put(130,170){\line(1,0){10}}
\put(140,170){\circle*{5}}
\put(138,175){$l$}

\put(210,170){$B_l$}
\put(240,170){\circle*{5}}
\put(238,175){1}
\put(240,168){\line(1,0){20}}
\put(240,171){\line(1,0){20}}
\put(246,167){$<$}
\put(260,170){\circle*{5}}
\put(258,175){2}
\put(260,170){\line(1,0){30}}
\put(280,170){\circle*{5}}
\put(278,175){3}
\put(300,170){\circle*{1}}
\put(310,170){\circle*{1}}
\put(320,170){\circle*{1}}
\put(330,170){\line(1,0){10}}
\put(340,170){\circle*{5}}
\put(338,175){$l$}
\end{picture}
\setlength{\unitlength}{\tmplen}
}
\end{center}
\caption{\label{tab:dynkin} Dynkin diagram of irreducible root systems}

\end{table}

\subsection{A parameterization of the irreducible representations in
defining characteristic}\label{subsec:param}

Let  $\Galg$,  $\Galg_{sc}$  and  $K$  be  as  above.  We  want  to  give  a
parameterization of the  irreducible defining characteristic representations
of  $\Galg^F$  by   describing  the  parameter  sets~(A),   (B)  and~(C)  of
theorem~\ref{maintheorem}.

Since $\Galg$ is semisimple we have $\Galg = \Galg'$ and so the group~(C) is
trivial in all cases considered here. 

The parameter sets~(A) and~(B) only depend on $K^F$,  the $F$-fixed elements
of the kernel of the isogeny $\Galg_{sc} \to \Galg$.

We now consider the possibilities for $\Galg$, $K$ and $K^F$ for the various
types of root systems separately.

We  make  use  of  the result  in~\cite[\S6.2]{luebeckdef}  which  describes
explicitly the elements  of the center $Z$ of $\Galg_{sc}$  in all cases (as
elements of $(\Q_{p'}/\Z)^l$ as explained above).

For any positive integers $l$ and  $q$, we will write $\mathcal E_{l,q}$ for
the  set of  tuples $(\lambda_1,\ldots,\lambda_l)\in\Z^l$  such that  $0\leq
\lambda_i<q$ for all $1\leq i\leq l$.

\subsubsection{Type $A_l$}
The group $Z$ is cyclic of order $m=(l+1)_{p'}$, generated by
\[z  =   \left(\frac{1}{m},  \frac{2}{m},  \ldots,   \frac{l}{m}\right)  \in
(\Q_{p'}/\Z)^l.\]
For each divisor $e$ of $l+1$ there  is an algebraic group $\Galg$ such that
the index  of $\Z R  \leq X$ is  $e$, we denote  its type by  $(A_l)_e$. So,
$e=l+1$  yields the  simply connected  groups, isomorphic to $\SL_{l+1}(\bar
k)$, and  $e=1$ yields  the adjoint  groups, isomorphic  to $\PGL_{l+1}(\bar
k)$.

Assume that  $\Galg$  is of type $(A_l)_e$.

Then $K$  is the subgroup of  $Z$ of order $((l+1)/e)_{p'}  = m/e_{p'}$ (the
group generated by $e_{p'}z$).

For the  Frobenius morphism  $F_\epsilon$ on $\Galg_{sc}$  and $i\in  \Z$ we
have $F_\epsilon(iz)  = iz$ if  and only if  $(q-\epsilon)i \in m\Z$  if and
only if $(m/\gcd(m,q-\epsilon)) \mid i$.

Combining, we find that the group $K^{F_\epsilon}$ is the subgroup of $Z$ of
order 
\[d := \gcd(m/e_{p'}, \gcd(m, q-\epsilon)) = \gcd(m/e_{p'}, q-\epsilon)
= \gcd((l+1)/e,  q-\epsilon)\] 
(the  last equation  because $  q-\epsilon$ is prime to $p$).

We   have  found   that  the   parameter  group~(B)   is  cyclic   of  order
$d$.  The  set~(A)   consists  of  the  $q$-restricted   weights  which  are
trivial  on  the  generator  $(m/d)z$ of  $K^{F_\epsilon}$.  These  are  the
$(\lambda_1,\ldots,\lambda_l)\in\mathcal E_{l,q}$ such that
\begin{equation}
\label{eq:Al}
\sum_{i=1}^{l} i \lambda_i\equiv 0\mod d.
\end{equation}

\subsubsection{Types $B_l$ and $C_l$}
In these two cases, the center $Z$ of $\Galg_{sc}$ has order $m=\gcd(2,p+1)$
and has the following generators.
\[\renewcommand{\arraystretch}{1.4}\begin{array}{c|c}
\textrm{Type}&\textrm{Generator}\\
\hline
B_l&(\frac{1}{m}, 0,\ldots,0)\\
C_l&(\frac{l}{m},\frac{l-1}{m},\ldots,\frac{2}{m},0,\frac{1}{m})
\end{array}\]\vspace{2pt}

There are two possibilities for $\Galg$, the simply-connected type where $K$
and so  $K^F$ are trivial, and  the adjoint type  where $K = Z$  and clearly
$K^F=K$ (since $K$ is of order $1$ or $2$).

So, when $p=2$  or $\Galg$ is simply-connected then  the parameter group~(B)
is  trivial  and  the  parameter  set~(A)  consists  of  all  $q$-restricted
weights. Otherwise, for  odd $q$ and $\Galg$ of adjoint  type, the group~(B)
is  of  order~$2$  and  the   parameter  set~(A)  consists  of  the  weights
$(\lambda_1,\ldots,\lambda_l)\in\mathcal  E_{l,q}$ satisfying  the following
equation. In the case of type $B_l$, the equation is
\begin{equation}
\label{eq:Bl}
\lambda_1\equiv 0\mod 2,
\end{equation} 
and in case of type $C_l$ it depends on the parity of $l$. This is
\begin{equation}
\label{eq:Cl}
\sum_{1\leq i\leq l,\,i\ \textrm{even}}\lambda_i\equiv 0\mod
2\qquad\textrm{or}\quad
\sum_{1\leq i\leq l,\,i\ \textrm{odd}}\lambda_i\equiv 0\mod
2,
\end{equation}
according to $l$ being odd or even. 

\subsubsection{Type $D_l$, $l \geq 4$}

Assume that $l=2k+1$  is odd. Then the  center $Z$ is cyclic  of order $m=4$
for odd $p$ and $m=1$ for $p=2$, and is generated by
\[z=\left(\frac{1}{m},     \frac{3}{m},    \frac{2}{m},0,     \frac{2}{m},0,
\frac{2}{m}, \ldots, 0, \frac{2}{m}\right).\]

There are  three possibilities  for $\Galg$,  the simply-connected  type for
which $K$  is trivial, or  $\Z R$ is  of index $2$  in $X$, then  $\Galg$ is
isomorphic to $\operatorname{SO}_{2l}(\bar k)$ and $K$ is generated by $2z$,
or the group of adjoint type where $K = Z$.

So, if $p=2$ or if $\Galg$ is simply-connected then the parameter group~(B)
is trivial and~(A) consists of all $q$-restricted weights.

If   $p$   is   odd   and  $\Galg$   of   type   $\operatorname{SO}$,   then
$K^{F_\epsilon}  =  K$  (since  $2z$  is  the  only  element  of  order  $2$
in  $Z$),  so the  group~(B)  is  of order  $2$.  In  this case, $\Galg^F  =
\operatorname{SO}_{2l}^{\epsilon}(q)$,  the  parameter  set~(A)  consists of
the weights $(\lambda_1,\ldots,\lambda_l)\in
\mathcal E_{l,q}$ such that
\begin{equation}
\label{eq:SOimpair}
\lambda_1+\lambda_2\equiv 0 \mod 2.
\end{equation}

Let $p$  be odd  and $\Galg$  be of  adjoint type,  then $K=Z$.  If $q\equiv
\epsilon\mod 4$  then $K^{F_\epsilon}  = K$, so  the parameter  group~(B) is
cyclic  of order  $4$  and the  parameter set~(A)  consists  of the  weights
$(\lambda_1,\ldots,\lambda_l)\in \mathcal E_{l,q}$ such that
\begin{equation}
\label{eq:Dlimpair}
\sum_{i=1}^k2\lambda_{2i+1}\equiv \lambda_2-\lambda_1\mod 4.
\end{equation}
Otherwise,  if  $q\equiv  -\epsilon\mod   4$  then  $K^{F_\epsilon}$  is  of
order~$2$  and  the parameter  sets~(B)  and~(A)  are  the  same as  in  the
$\operatorname{SO}$-case.

\subsubsection{Type $D_l$, $l \geq 4$}
Assume now  that $l=2k$  is even.  Then $Z$ is  elementary abelian  of order
$4$ if $p$ is odd and  trivial if $p=2$. 
If $p$ is odd then $Z$ is generated by
\[z_1 =  \left(\frac{1}{2}, 0,  0, \frac{1}{2},  0, \frac{1}{2},  \ldots, 0,
\frac{1}{2}\right)\quad\textrm{ and  }\quad z_2  = \left(0,  \frac{1}{2}, 0,
\frac{1}{2}, 0, \frac{1}{2}, \ldots, 0, \frac{1}{2}\right),\]
if $p=2$ we set $z_1=z_2=1$.

Here, for any $q$, $F_1$ is the identity on $Z$, $F_{-1}$ permutes $z_1$ and
$z_2$,  and in  case $l=4$  the Frobenius  $F_3$ permutes  $z_1$, $z_2$  and
$z_1+z_2$ cyclically.

If $\Galg$ is simply-connected  or $p=2$, then $K = K^F  = 1$, the parameter
group~(B) is trivial and the set~(A) consists of all $q$-restricted weights.

If  the index  of $\Z  R$ in  $X$ is  $2$, there  are two  possibilities for
$\Galg$. Either  $\Galg$ has  only Frobenius morphisms  of type  $F_1$, then
$\Galg$ is  isomorphic to a half  spin group $\operatorname{HSpin}_{2l}(\bar
k)$ and $K = K^F$ is generated by $z_1$ (or by $z_2$). In this case, for odd
$p$, the parameter group~(B) is of order $2$ and the set~(A) consists of the
weights $(\lambda_1,\ldots,\lambda_l)\in\mathcal E_{l,q}$ such that
\begin{equation}
\label{eq:HS}
\lambda_2+\lambda_4+\cdots+\lambda_{2k}\equiv 0\mod 2.
\end{equation}
Otherwise,  $K$  is  generated   by  $z_1+z_2$  and  $K=K^{F_\epsilon}$  for
$\epsilon  \in   \{\pm1\}$.  Then  $\Galg$   is  isomorphic  to   a  special
orthogonal group  $\operatorname{SO}_{2l}(\bar k)$  and $\Galg^{F_\epsilon}$
is  isomorphic  to  $\operatorname{SO}_{2l}^\epsilon(q)$. For  odd  $p$  the
parameter group~(B)  is also of  order $2$ and  the set~(A) consists  of the
weights $(\lambda_1,\ldots,\lambda_l)\in\mathcal E_{l,q}$ such that
\begin{equation}
\label{eq:SOpair}
\lambda_1+\lambda_2\equiv 0\mod 2.
\end{equation}

The final possibility is  that $\Galg$ is of adjoint type and  $K = Z$. Then
$K^{F_{-1}}$ is generated by $z_1+z_2$ and  for odd $p$ and $\epsilon=-1$ we
get the same parameter sets~(B)  and~(A) as in the $\operatorname{SO}$-case.
Furthermore, we  have $K^{F_1}  = K$,  so for odd  $p$ and  $\epsilon=1$ the
parameter group~(B)  is elementary  abelian of  order~$4$ and  the parameter
set~(A)  consists  of the  weights  $(\lambda_1,\ldots,\lambda_l)\in\mathcal
E_{l,q}$ such that
\begin{equation}
\label{eq:Dlpair}
\left\{
\renewcommand{\arraystretch}{1.4}
\begin{array}{l}
\displaystyle{\sum_{i=2}^k\lambda_{2i}\equiv \lambda_1\mod 2,}\\
\lambda_1=\lambda_2\mod 2.
\end{array}\right.
\end{equation}
If $l=4$ then $K^{F_3}  = 1$ and we get the same  parameterization as in the
simply connected case for $F=F_3$.

\subsubsection{Types $G_2$, $F_4$ and $E_8$}
In these cases $Z$ and so $K=K^F$ and the parameter group~(B) are trivial.
The set~(A) consists of all $q$-restricted weights.

\subsubsection{Type $E_6$}
The group $Z$ is cyclic of order $m=3$  if $p\neq 3$ and $m=1$ if $p=3$,  it
is  generated by  $z=(\frac{1}{m},0,\frac{2}{m},0,\frac{1}{m},\frac{2}{m})$.
The group $\Galg$ can either be simply-connected or of adjoint type.

If  $\Galg$  is simply-connected  or  $p=3$  then $K=K^{F_\epsilon}=1$,  the
parameter  group~(B)  is  trivial  and~(A) consists  of  all  $q$-restricted
weights.

If $\Galg$ is of adjoint type then $K=Z$. We have $K^{F_\epsilon} = K$ if $q
\equiv \epsilon \mod  3$. In that case the parameter  group~(B) is cyclic of
order~$3$ and the parameter set~(A) consists of the weights
$(\lambda_1,\ldots,\lambda_6)\in\mathcal E_{6,q}$ such that
\begin{equation}
\label{eq:E6}
\lambda_1-\lambda_3+\lambda_5-\lambda_6\equiv 0\mod 3.
\end{equation}
For  $q  \equiv -\epsilon  \mod  3$  we  have  $K^{F_\epsilon} =1$  and  the
parameter sets are as in the simply-connected case.

\subsubsection{Type $E_7$}
The group $Z$ is  cyclic of order $m=2$ if $p\neq 2$ and  $m=1$ if $p=2$, it
is generated by $z=(0,\frac{1}{m},0,0,\frac{1}{m},0,\frac{1}{m})$. The group
$\Galg$ is either simply-connected or of adjoint type.

If  $\Galg$  is  simply-connected  or  if $p=2$  then  $K=K^F$  is  trivial,
the  parameter  group~(B)  is  trivial  and  the  set~(A)  consists  of  all
$q$-restricted weights.

If  $\Galg$ is  of  adjoint  type then  $K  =  K^F =  Z$.  For  odd $p$  the
parameter group~(B) is of order~$2$ and  the set~(A) consists of the weights
$(\lambda_1,\ldots,\lambda_7)\in\mathcal E_{7,q}$ such that
\begin{equation}
\label{eq:E7}
\lambda_2+\lambda_5+\lambda_7\equiv 0\mod 2.
\end{equation}

\subsection{Application: number of semisimple classes}

In  this section,  we  will compute  the number  of  isomorphism classes  of
irreducible $\bar{\F}_p$-modules (or, equivalently, the number of semisimple
classes)  of the  finite groups  $\Galg^F$ for  all simple  algebraic groups
$\Galg$ defined over $\F_q$.

\begin{theorem}
\begin{itemize}
\item[(a)]  Let  $\Galg$  be  a  connected  reductive  group  of  semisimple
rank  $l$, such  that its  derived group  $\Galg'$ is  simply-connected. Let
$Z(\Galg)^0$ be the connected component of  the center of $\Galg$. We assume
that  $\Galg$ is  defined  over  $\F_q$ and  denote  $F: \Galg\to\Galg$  the
corresponding Frobenius  morphism. Then  the number of  semisimple conjugacy
classes of $\Galg^F$ is $q^l |(Z(\Galg)^0)^F|$. In particular, if $\Galg$ is
simply-connected this number is $q^l$.
\item[(b)] Now  let $\Galg$ be  a simple  connected reductive group  of rank
$l$, defined over $\F_q$ with corresponding Frobenius morphism $F$. Then the
number of semisimple  conjugacy classes of $\Galg^F$ is either  $q^l$, or it
is given in table~\ref{tab:semisimple}.
\end{itemize}
\end{theorem}

\textbf{Proof.}
(a) This follows from theorem~\ref{maintheorem}. Under the given assumptions
the   parameter  set~(B)   is  trivial,   and  the   set~(A)  contains   all
$q$-restricted weights. The parameter set~(C) contains $|(\Galg/\Galg')^F| =
|(Z(\Galg)^0)^F|$ elements. See also~\cite[3.7.6(ii)]{Ca85} for a completely
different proof of this result.

(b) This  will be  shown in  the rest  of this  section. Here  the parameter
group~(C)  is  always  trivial.  We  need   to  go  through  all  the  cases
of  subsection~\ref{subsec:param}.   Whenever  the  group~(B)   is  trivial,
the  set~(A)  consists of  the  $q^l$  elements in  $\mathcal{E}_{l,q}$.  In
table~\ref{tab:semisimple}  we collect  the cases  with non-trivial~(B)  and
find the cardinalities of the sets~(A)  by counting the solutions of certain
modular equations.
\hfill$\Box$

We  denote   by  $\Lambda$   the  set  of   parameters~(A)  for   the  group
$\Galg^F=(\Galg_{\operatorname{sc}}/K)^F$. By theorem~\ref{maintheorem}, the
number  of  isomorphism  classes   of  irreducible  $\bar{\F}_p$-modules  of
$\Galg^F$  is $|K^F|\cdot|\Lambda|.$  

The following lemma will be useful in several cases.

\begin{lemma}
Assume  that   $q$  is  odd,   and  for   any  positive  integers   $n$  and
$\nu\in\{0,1\}$,  define  $$E_{n,\nu}=\left\{(\lambda_1,\ldots,\lambda_n)\in
\Z^n\,|\,0\leq \lambda_i\leq q-1,\quad \sum_{i=1}^{n}\lambda_i\equiv \nu\mod
2\right\}.$$ Then, we have
\[|E_{n,\nu}|=\frac{q^n+1-2\nu}{2}\]
\label{nbpair}
\end{lemma}
\textbf{Proof.}
This follows easily by induction on $n$.
\hfill$\Box$

\begin{table}
\[\renewcommand{\arraystretch}{1.4}\begin{array}{c|c|c|c|c|c}
\textrm{Type}&K^F&F&\Galg&\textrm{condition}&
|\textrm{semisimple classes}|\\
\hline
A_l&\Z_{d}&F_\epsilon&(A_l)_e&d=\gcd(\frac{l+1}{e},q-\epsilon)&
\sum_{d' \mid d} \varphi(d') q^{(l+1)/d'-1}\\
B_l&\Z_2&&\textrm{adjoint}&p\neq 2&q^l+q^{l-1}\\
C_l&\Z_2&&\textrm{adjoint}&p\neq 2&q^l+q^{\lfloor l/2\rfloor}\\
D_{l},l\textrm{ even}&\Z_2^2&F_1&\textrm{adjoint}&p\neq 2&q^{l}+q^{l-2}
+2q^{l/2}\\
&\Z_2&F_{-1}&\textrm{adjoint}&p\neq 2&q^{l}+q^{l-2}\\
&\Z_2&F_{\epsilon}&\operatorname{SO}
&p\neq 2&q^{l}+q^{l-2}\\
&\Z_2&F_{1}&\operatorname{HSpin}&p\neq 2&q^{l}+q^{l/2}\\
D_{l}, l\textrm{ odd}&\Z_4&F_{\epsilon}&\textrm{adjoint}
&q\equiv \epsilon \!\!\!\mod 4&
q^{l}+q^{l-2}+2q^{(l-3)/2}\\
&\Z_2&F_{\epsilon}&\textrm{adjoint}
&q\equiv -\epsilon \!\!\!\mod 4&
q^{l}+q^{l-2}\\
&\Z_2&F_{\epsilon}&\operatorname{SO}&p\neq 2&
q^{l}+q^{l-2}\\
E_6&\Z_3&F_{\epsilon}&\textrm{adjoint}
&q\equiv \epsilon \!\!\!\mod  3&q^6+2q^2\\
E_7&\Z_2&&\textrm{adjoint}&p\neq 2&q^7+q^4\\
\end{array}
\]
\caption{\label{tab:semisimple}Number of semisimple classes ($\varphi$ is
the Euler $\varphi$-function)}
\end{table}

\textit{Types $B_l$ and $C_l$.}\quad
We must  consider the  case that  $p\neq 2$  and $K=Z$.  Then $\Galg$  is of
adjoint type. If $\Galg$ is of  type $B_l$ we use equation~(\ref{eq:Bl}) and
obtain
\[\Lambda=\{(\lambda_1,\ldots,\lambda_l)\in \mathcal
E_{l,q}\,|\,\lambda_l\in 2\Z\}.\]
Thus,  lemma~\ref{nbpair} gives  $|\Lambda|=q^{l-1}\cdot\frac{q+1}{2}$. Now,
since  $|K^F|=2$, the  entry  for case  $B_l$ in  table~\ref{tab:semisimple}
follows.

If  $\Galg$  is  of  type  $C_l$  with  $l=2k$  (resp.  $l=2k+1$),  then  in
equation~(\ref{eq:Cl}) there are $k$ summands  (resp. $k+1$ summands) in the
sum. Hence, lemma~\ref{nbpair} gives
\[|\Lambda|=q^k\cdot\frac{q^k+1}{2}\quad\left(\textrm{resp. }
q^k\cdot\frac{q^{k+1}+1}{2}\right).\]
Since  $k=\lfloor l/2\rfloor$  and $|K^F|=2$,  the entry  for type  $C_l$ in
table~\ref{tab:semisimple} follows.

\textit{Type   $D_l$.}\quad   We   only    need   to   consider   the   case
$p   \neq  2$.   First   assume   that  $l=2k$.   We   compute  the   number
of   elements   in   the   set~(A)   for   $\Galg_{\operatorname{ad}}^{F_1}$
using    equation~(\ref{eq:Dlpair})    as    follows.    If    $\lambda_{1}$
is    odd,    then    $\lambda_{2}$    is    odd.    This    implies    that
$\lambda_4+\lambda_6+\cdots+\lambda_{2k}\in  2\Z+1$. By  lemma~\ref{nbpair},
there    are   $q^{k-1}\left(\frac{q-1}{2}\right)^2\cdot\frac{q^{k-1}-1}{2}$
such           solutions.           Similarly,           there           are
$q^{k-1}\left(\frac{q+1}{2}\right)^2\cdot\frac{q^{k-1}+1}{2}$ solutions such
that $\lambda_{1}$ is even. Therefore, we deduce that
\[|\Lambda|=q^{k-1}\cdot\frac{q^{k+1}+2q+q^{k-1}}{4}.\]
In the same way, using lemma~\ref{nbpair},  we count the number of solutions
of  equations~(\ref{eq:SOpair})  and~(\ref{eq:HS}), giving  $|\Lambda|$  for
$\operatorname{SO}_{2l}^{\epsilon}(q)$  and  $\operatorname{HSpin}_{2l}(q)$,
respectively.

Suppose   now  that   $l=2k+1$.  Assume   that  $K=Z$   and  that   $q\equiv
\epsilon\mod   4$.  By   equation~(\ref{eq:Dlimpair})   we   have  to   find
the   number    of   solutions   $(\lambda_1,\ldots,\lambda_{2k})\in\mathcal
E_{2k,q}$   of    $2(x_3+x_5+\cdots+x_{2k+1})+x_{1}-x_{2}\in   4\Z$.   There
are  $q^{k-1}.\frac{q^k+1}{2}$   tuples  $(\lambda_1,\ldots,\lambda_{2k-1})$
with        $0\leq        \lambda_i\leq        q-1$,        such        that
$(\lambda_3+\lambda_5+\cdots+\lambda_{2k+1})$  is even.  For each  tuple, we
have to  find the number  of solutions of  $\lambda_{1}-\lambda_{2}\in 4\Z$.
There  are $\left(\frac{q+3}{4}\right)^2+3\cdot\left(\frac{q-1}{4}\right)^2$
such            solutions.             Thus,            there            are
$n_0=\frac{1}{8}\cdot       q^{k-1}\cdot(q^{k+2}+3q^k+q^2+3)$      solutions
$(\lambda_1,\ldots,\lambda_{2k+1})\in\mathcal    E_{2k+1,q}$,   such    that
$\lambda_3+\lambda_5\cdots+\lambda_{2k+1}$   is   even.   Similarly,   there
are    $n_1=\frac{1}{8}\cdot   q^{k-1}\cdot(q^{k+2}-q^k-q^2+1)$    solutions
$(\lambda_1,\ldots,\lambda_{2k+1})\in\mathcal    E_{2k+1,q}$,   such    that
$\lambda_1+\lambda_3\ldots+\lambda_{2k-1}$  is  odd.   Adding up we find for
the case $|K^{F_\epsilon}| = 4$ that
\[|\Lambda|=n_0+n_1=\frac{q^{k-1}}{4} (q^{k+2}+q^k+2).\]

Finally,  again   with  lemma~\ref{nbpair},   we  count  the   solutions  of
equation~(\ref{eq:SOimpair})   and   obtain   $|\Lambda|$  for   the   cases
$\operatorname{SO}_{2l}^{\epsilon}(q)$,  and for  the cases with $\Galg$  of
adjoint type and $q\equiv -\epsilon\mod 4$.

\textit{Types $E_6$ and $E_7$.}\quad
We   only   need   to   consider   $\Galg$  of   adjoint   type.   In   type
$E_6$  with  $q\equiv   \epsilon  \mod  3$  we   compute  $|\Lambda|$  using
equation~(\ref{eq:E6})  and  lemma~\ref{lemma:count}.  For  type  $E_7$  and
odd   $p$   we   conclude   as  above   using   equation~(\ref{eq:E7})   and
lemma~\ref{nbpair}.

\textit{Type $A_l$.}\quad  We have postponed this  case because it is  a bit
trickier  to derive  a closed  formula for  the cardinality  $|\Lambda|$. We
need  to count  the solutions  of equation~(\ref{eq:Al})  to find  the first
line  of table~\ref{tab:semisimple}.  Instead of  counting the  solutions of
equation~(\ref{eq:Al})  we  can  introduce another  coordinate  $\lambda_0$,
count the solutions of the equation
\[ \sum_{i=0}^l i\lambda_i \equiv 0 \mod d \]
with $0\leq \lambda_i  < q$ for $0\leq  i \leq l$, and divide  the result by
$q$.

The number of solutions of this modified equation follows from the following
lemma applied with $n=l+1$ and $m = d$.

\begin{lemma}\label{lemma:an}
Let  $n\geq  2$  be  an   integer,  and  $m|n$.  Let  $\epsilon\in\{-1,1\}$.
For    any     positive    integer    $q$     with    $m    \mid     (q    -
\epsilon)$   define   $E=\{0,\ldots,q-1\}$    and   write   $X=\{(\lambda_0,
\lambda_1,\ldots,\lambda_{n-1})\in     E^{n}\,|\,    \sum\limits_{i=0}^{n-1}
ix_i\equiv 0\mod m\}$. Then
\[|X|=\frac{1}{m}\sum_{d|m}\varphi\left(d\right)q^{n/d}.\]
\end{lemma}

\textbf{Proof. }
The following proof was shown to us by Darij Grinberg, who generously
allowed us to include it in this article.

For a non-negative integer $t$ let 
\[X_t = \{(\lambda_0,\ldots,
\lambda_{n-1}) \in E^n\mid\; \sum_{i=0}^{n-1} i \lambda_i = t\}.\]
We want to investigate $|X| = \sum_{t\in m\Z} |X_t|$. 

For integers $j$ with $0\leq j\leq n-1$ we consider the polynomials
\[P_j(z) = 1 + z^j +z^{2j} + \ldots +z^{(q-1)j} \in \mathbb{C}[z].\]
It is clear that $|X_t|$ is the coefficient of $z^t$ in the product of 
the $P_j(z)$:
\[P(z) = \prod_{j=0}^{n-1} P_j(z) = \sum_{t=0}^\infty |X_t|z^t.\]

Now   let  $\zeta   \in  \mathbb{C}$   be   a  primitive   $m$-th  root   of
unity.  We  will  use  repeatedly  the   fact  that  for  $t\in\Z$  the  sum
$\sum_{k=0}^{m-1}\zeta^{tk}$  equals  $m$  if   $m\mid  t$  and  equals  $0$
otherwise (use the formula for geometric sums).

We evaluate
\[\sum_{k=0}^{m-1} P(\zeta^k) = 
\sum_{k=0}^{m-1} \sum_{t=0}^\infty |X_t|\zeta^{k t} =
\sum_{t=0}^\infty |X_t| \sum_{k=0}^{m-1} (\zeta^{t})^k =
\sum_{\substack{t=0\\ m\mid t}}^\infty |X_t| \cdot m =
m |X|.\]

From now we fix a $k\in\Z$ with $0\leq k \leq m-1$ and set
$d = \frac{m}{\gcd(m,k)}$.

We will show that 
\[ P(\zeta^k) = q^{n/d}.\]

This proves the lemma, because for $d\mid m$ we have
\[|\{0\leq k < m\mid\; d = m/\gcd(m,k)\}| = |\frac{m}{d}\cdot \{
0\leq i < d\mid\; \gcd(i,d) = 1\}| = \varphi(d).\]

It is easy to evaluate each $P_j(z)$ at $\zeta^k$ for 
$0\leq j \leq n-1$:
\[P_j(\zeta^k) = 1 + \zeta^{jk} + (\zeta^{jk})^2 + \ldots +
(\zeta^{jk})^{q-1} = \left\{ \begin{array}{ll}
q,& \textrm{if } m\mid jk\\
1,& \textrm{if } m\nmid jk \textrm{ and } \epsilon=1\\
-\zeta^{-kj},& \textrm{if } m\nmid jk \textrm{ and } \epsilon=-1
\end{array}\right. 
\]
because for  $m\nmid jk$ every  $m$ consecutive summands  sum up to  $0$. In
case $m\mid  (q-1)$ it remains  the last summand which  is $1$. And  in case
$m\mid (q+1)$  an additional summand $(\zeta^{jk})^q  = \zeta^{-jk}$ cancels
all the previous ones.

For an integer $k$, $d = m/\gcd(m,k)$ and $j\in\Z$ we have that 
$m\mid jk$ if and only if $d\mid j$. 

Since $m\mid n$ we also have $d\mid n$ and so there are $n/d$ indices
$j$ with $0\leq j \leq n-1$ and $P_j(\zeta^k) = q$. 

In case $m\mid (q-1)$ we have $P_j(\zeta^k) = 1$ for the remaining $j$ with
$d\nmid j$. Taking the product we get
\[ P(\zeta^k) = q^{n/d}. \]

To see that the same is true in case $m\mid (q+1)$ we must show that
\[\prod_{\substack{j=0\\ d\nmid j}}^{n-1} (-\zeta^{-jk}) = 1.\]

The root of unity $\zeta^k$ and so also $\zeta^{-k}$ has order
$d$ ($= m/\gcd(m,k)$). So, if $d\mid j$ we have $(\zeta^{-k})^j = 1$ and
we get
\[\prod_{\substack{j=0\\ d\nmid j}}^{n-1} (-\zeta^{-jk}) = (-1)^{n-n/d}
\cdot \prod_{j=0}^{n-1} (\zeta^{-k})^j = (-1)^{n-n/d} \cdot
(\zeta^{-k})^{n(n-1)/2}.\]
Since $d\mid n$, we have: $d\nmid n(n-1)/2$  iff ($d$ is even and $(n/d)$ is
odd)  iff $(n/d)(d-1)  =  n -  n/d$  is  odd, and  in  this case  $(d/2)\mid
n(n-1)/2$  so   that  $(\zeta^{-k})^{n(n-1)/2} =  -1$. This  shows that  the
right hand side in the last displayed equation is always $1$.

This proves the lemma.
\hfill$\Box$

\begin{remark}\rm
The  results given  in  table~\ref{tab:semisimple} are  not  new. They  were
worked  out in~\cite{BrGGC}  using sophisticated  results from  the ordinary
representation theory  of the groups  $\Galg^F$ in good  characteristic. The
completely different  approach in  this section is  more elementary  (and it
works for arbitrary root data and characteristics).

For  small  rank  groups,  in particular  the  exceptional  types,  detailed
parameterizations of all  conjugacy classes were computed,  this also yields
the number of semisimple conjugacy classes, see~\cite{LWeb}.
\end{remark}

%%  \bibliographystyle{plain}
%%  \bibliography{references}

\end{document}